\documentclass[11pt,reqno]{amsart}
\usepackage{amsfonts,amsmath,amsthm,amssymb,bm}
\usepackage[usenames,dvipsnames]{color}
\usepackage{enumerate,nicefrac}
\usepackage{hyperref,multicol}
\usepackage{rotating}
\usepackage[labelsep=space]{caption}
\usepackage[all]{xy}
\usepackage{tikz}
\usetikzlibrary{decorations.markings}
\usepackage[labelformat=empty]{caption}
\usepackage{mathtools}

\newtheorem{theorem}{Theorem}[section]
\newtheorem{lemma}[theorem]{Lemma}

\theoremstyle{definition}

\usepackage{multirow, array, graphicx}
\usepackage{float}

\usepackage[ansinew]{inputenc}
\usepackage{ifthen}
\usepackage{animate}

\newcommand{\be}{\begin{equation}}
\newcommand{\ee}{\end{equation}}

\title[On a Variant of Pillai{'}s problem]{On a Variant of Pillai{'}s problem involving convergent denominators of quadratic irrationals}

\author{Mohit Mittal}
\address{Department of Mathematics\\
Birla~Inst{i}tute~of~Technology~and~Science, Pilani 333\,031 \textsc{India}}
\email{mohit.mittal\symbol{64}pilani.bits-pilani.ac.in}

\begin{document}

\keywords{Diophantine equation, Pillai's  problem, linear recurrence relation, Baker's method, multiplicatively dependent numbers, reduction method}

\subjclass[2020]{11D45, 11D61, 	11J86, 11R04, 11R11}

\begin{abstract}
Let $(q_{\alpha, n})_{n \geq 0}$ be the sequence of convergent denominators to the simple continued fraction expansion of $\alpha$. For certain specific choices of $\alpha$, this sequence is a Lehmer sequence. In this paper, we show that there are only finitely many integers $c$ such that the equation $q_{\alpha, n} - q_{\beta, m} = c$ has at least two distinct solutions $(n,m)$, where $\alpha,\beta$ are quadratic irrationals with $\mathbb{Q}(\alpha)\neq \mathbb{Q}(\beta)$. In specific instances, we solve the equation $q_{\alpha, n} - q_{\beta, m} = c$ completely and explicitly list all solutions.
\end{abstract}

\maketitle

\section{Introduction}

\indent \ \ In 1931, Pillai \cite{Pillai_31} proved that for fixed positive integers $a$ and $b$, with $\min(a, b) \geq 2$ and $\nicefrac{\log a}{\log b}$ irrational, the number of solutions $(x,y)$ to the Diophantine inequalities $0 <  a^{x} - b^{y} \leq c$ is asymptotically equal to \begin{equation*}
    \frac{(\log c)^{2}}{2 (\log a)(\log b)} 
\end{equation*} as $c \to \infty$. This work drew Herschfeld's \cite{HF_35, HF_36} attention, who proved, in 1935, that the Diophantine equation \begin{equation*}
2^{x} - 3^{y} = c, \end{equation*} has at most one solution when $|c|$ is sufficiently large and fixed.
\noindent Following this, in 1936, Pillai \cite{Pillai_36, Pillai_37} considered the more general Diophantine equation \begin{equation}\label{Geneqn_Pillai}
    a^{x} - b^{y}=c,
\end{equation} where $a, b, c$ are fixed nonzero integers with $a>b \geq 1$. He showed that for coprime $a$ and $b$, sufficiently large values of $c$ allow at most one positive integer solution $(x,y)$. Further, he put forth two conjectures. The first conjecture states that the only integers $c$ that admit at least two representations of the form $2^{x} - 3^{y}$ come from the set $\left\{-1, 5, 13\right\}$, with explicit representations 
\begin{align*}
  2^{3}-3^{2}   &= -1 = 2-3,\\
  2^{5} - 3^{3} &= 5  = 2^{3} - 3,\\
  2^{8} - 3^{5} &= 13 = 2^{4} - 3.
\end{align*} This conjecture was confirmed by Stroeker and Tijdeman \cite{Stroeker_Tijdeman_82} in 1982. (See also \cite{Bennett_01}.)

Secondly, he conjectured that, for fixed nonzero $c$, equation \eqref{Geneqn_Pillai} has at most finitely many solutions in integers $a, b, x,$ and $y$, all greater than one. The case $c=1$ is  Catalan's conjecture, one of the famous classical problems in number theory. In 1976, Tijdeman \cite{Tijdeman_76} proved that equation \eqref{Geneqn_Pillai} with $c=1$ has only finitely many solutions. Catalan's conjecture was proved in 2002 by Mih\u{a}ilescu \cite{Mihail_04}. However, Pillai's conjecture remains unproven for $c \geq 2$. 

\indent In 2016, Chim, Pink \& Ziegler \cite{Chim_Ziegler_17} considered the Diophantine equation  \begin{equation}\label{CPZ_eqn1}
    F_{n} - T_{m} = c,
\end{equation}
where $c$ is a fixed integer and $(F_n)_{n \geq 0}$, $(T_m)_{m \geq 0}$ are the sequences of Fibonacci numbers and Tribonacci numbers, respectively. This type of equation can be viewed as a variation of Pillai's equation. 
They showed that the only integers $c$ having at least two distinct representations of the form $F_{n} - T_{m}$ come from the set $$\mathcal{C}=\left\{0,1,-1,-2,-3,4,-5,6,8,-10,11,-11,-22,-23,-41,-60,-271\right\},$$ and also gave all representations of each $c \in \mathcal{C}$.

In the same year, Chim, Pink \& Ziegler \cite{Chim_Ziegler_18} considered the Diophantine equation 
\begin{equation}\label{CPZ_eqn2}
    U_{n} - V_{m} = c,
\end{equation} where $c$ is a fixed integer and $(U_n)_{n \geq 0}$, $(V_m)_{m \geq 0}$ are linear recurrence sequences satisfying the \textit{dominant root condition}. By finding an effective upper bound for $|c|$, they showed that there are finitely many choices of $c$ which have at least two distinct representations of the form $U_{n} - V_{m}$. We refer to \cite{Bravo_Luca_Yazan_017, Dama_019_P, Dama_019,  Dama_020, Dama_Luca_Gomez_017, Dama_Luca_020, Dama_Luca_Rako_017, Lomeli_Hernandez_019, Hernandez_Luca_Rivera_019, Luca_Rihane_Hernane_019} and the references therein to get more insight into the theme of this paper.

In this paper, for a fixed integer $c$, we study the equation \begin{equation*}
    q_{\alpha, n} - q_{\beta, m} = c,
\end{equation*} which is a variation of Pillai's equation. Here, $( q_{\alpha, n} )_{n \geq 0}$ and $( q_{\beta, m} )_{m \geq 0}$ are the sequences of convergent denominators arising from the continued fraction expansions of quadratic irrationals $\alpha$ and $\beta$, respectively. This problem has not been addressed previously, as the linear recurrence sequences do not satisfy the dominant root condition (see Remark \ref{Remark}).

For specific choices of quadratic irrationals $\alpha$ and $\beta$, the sequences $( q_{\alpha, n} )_{n \geq 0}$ and $( q_{\beta, m} )_{m \geq 0}$ become Lehmer sequences \cite{bala_14} of the type  
\[
L_n(R,-1) =
\begin{cases} 
\frac{\xi^n - \eta^n}{\xi - \eta} & \text{if } n \text{ is odd}, \\
\frac{\xi^n - \eta^n}{\xi^2 - \eta^2} & \text{if } n \text{ is even}.
\end{cases}
\]
Here, $R \in \mathbb{Z_+}$ and $\xi, \eta$ are the roots of the quadratic equation $x^2 - \sqrt{R}x-1=0$. In these cases, we explicitly solve certain Diophantine equations and exhibit all of their solutions. We begin with the following theorem:
\begin{theorem}\label{main thm}
   Let $\alpha,\beta$ be real quadratic irrational numbers with simple continued fraction expansions 
    \begin{align*}
        \alpha&=[a_{0,0};\ldots,a_{0,r_0-1},\overline{b_{0,0},\ldots,b_{0,s_0-1}}], r_0\geq 0, s_0\geq 1,\\
        \beta&=[a_{1,0};\ldots,a_{1,r_1-1},\overline{b_{1,0},\ldots,b_{1,s_1-1}}], r_1\geq 0, s_1\geq 1
    \end{align*}
    and $\mathbb{Q}(\alpha) \neq \mathbb{Q}(\beta)$. Let $(q_{\alpha,i})_{i\geq 0},(q_{\beta,i})_{i\geq 0}$ be the sequences of denominators of the convergents to the continued fraction expansions of $\alpha,\beta$, respectively. Then there exists an effectively computable finite set
    $\mathcal{C}$ such that the integer $c$ has at least two distinct representations of the form $q_{\alpha, N}-q_{\beta, M}$ if and only if $c \in \mathcal{C}$.
\end{theorem}

Let $a, b \in \mathbb{N}$ such that $2 \leq a < b$. For $\alpha = [0; \overline{1, a}]$ and $\beta = [0; \overline{1, b}]$, we explicitly compute an upper bound for the size of the solutions.
Moreover, for $2 \leq a < b \leq 5$, we effectively compute the set $\mathcal{C}$. We note that for $\alpha = [0; \overline{1, a}]$, the sequence of convergent denominators is a Lehmer sequence with $R=a$. We prove
\begin{theorem}\label{thm on bound}
Let $a, b \in \mathbb{N}$ with $ 2 \leq a < b$. Let $\alpha,\beta$ be real quadratic irrational numbers with simple continued fraction expansions 
    \begin{align*}
        \alpha =[0;\overline{1, a}],\ 
        \beta =[0;\overline{1,b}]
    \end{align*}
    and $\mathbb{Q}(\alpha) \neq \mathbb{Q}(\beta)$. Let $(q_{\alpha,i})_{i\geq 0},(q_{\beta,i})_{i\geq 0}$ be the sequences of denominators of the convergents to the continued fraction expansions of $\alpha,\beta$, respectively. If $c$ is an integer that has at least two distinct representations of the form $q_{\alpha, N}-q_{\beta, M}$, then
\[ \max (N, M) \leq 1.1793 \times 10^{46} \times (\log \theta_{\beta,1})^{6} \times (\log (9.06 \times 10^{14} \times (\log \theta_{\beta,1})^{2} ))^{3}.\]
\end{theorem}

In particular, for $a=2$ and $b=3$ we extend Theorem \ref{thm on bound} to find the effectively computable set $\mathcal{C}$. Therefore, we have the following theorem:
\begin{theorem}\label{thm explicit_C}
Let $\alpha = [0; \overline{1,2}]$ and $\beta=[0;\overline{1,3}]$. The only integers $c$ that have at least two distinct representations of the form $q_{\alpha, N} - q_{\beta, M}$ belong to the set $$\mathcal{C} = \left\{ -4,-1,0,10,37\right\}.$$
Furthermore, for each $c \in \mathcal{C}$, two distinct representations of the form $c=q_{\alpha, N} - q_{\beta, M}$ is shown below:
\begin{align*}       
  q_{\alpha, 5} - q_{\beta, 4} = 15-19 &=\boxed{-4}  = 1-5 =q_{\alpha,0} - q_{\beta,3},\\
 q_{\alpha, 3} - q_{\beta, 3} = 4-5 &=\boxed{-1}  = 3-4 =q_{\alpha,2} - q_{\beta,2},\\
 q_{\alpha, 3} - q_{\beta, 2} = 4-4 &=\boxed{0} = 1-1 =q_{\alpha,1} - q_{\beta,0},\\
q_{\alpha, 5} - q_{\beta, 3} = 15-5 &=\boxed{10}  = 11-1 =q_{\alpha,4} - q_{\beta,0},\\
q_{\alpha, 7} - q_{\beta, 4} = 56-19 &=\boxed{37}  = 41-4 =q_{\alpha,6} - q_{\beta,2}.
\end{align*}
\end{theorem}

\noindent For other choices of $\alpha = [0; \overline{1,a}]$ and $\beta=[0;\overline{1,b}]$, where \\ $(a, b) \in \left\{(2, 4), (2, 5), (3, 4), (3, 5), (4,5)\right\}$, the set $\mathcal{C}$ along with two distinct representations of each $c \in \mathcal{C}$ is shown in Section \ref{appendix}.

\textit{Remark 1.4.}\label{Remark} We note that for a real quadratic irrational number $\alpha$, the sequence of convergent denominators of its continued fraction expansion satisfies a linear recurrence relation, but the associated characteristic polynomial may not have a dominating root. For instance, if $\alpha=\sqrt{27}$, then its continued fraction expansion is $[5; \overline{5, 10}]$ and $q_{i+4} = 52q_{i+2}- q_{i}$, $i \geq 0$. The characteristic polynomial is $x^{4}-52x^{2}+1$ and its roots are $\pm \sqrt{26+15\sqrt{3}}, \pm \sqrt{26-15\sqrt{3}}$. To circumscribe this issue, we use an argument of Peth\H{o} \cite{Petho_81}
(see also Lenstra and Shallit \cite{LeSh93}) to split the sequence into finitely many subsequences, each satisfying the same binary recurrence relation (with different initial terms, see Section \ref{prelim}). \\
In Section \ref{prelim}, we mention some preliminary results. In Sections \ref{Sec_main thm}, \ref{Sec_thm on bound}, and \ref{Sec_thm explicit_C},  we prove Theorems \ref{main thm}, \ref{thm on bound}, and \ref{thm explicit_C}, respectively.

\section{Preliminaries}\label{prelim}
Let $\alpha$ be a real quadratic irrational number. Let its simple continued fraction expansion be $\alpha=[a_0;\ldots,a_{r-1},\overline{b_0,\ldots,b_{s-1}}]$, $r\geq 0, s\geq 1$.
Let $(q_{\alpha,i})_{i\geq 0}$ be the sequence of denominators of the convergents to $\alpha$. Then by Peth\H{o} \cite[Lemmata 1 \& 2]{Petho_81} (see also
\cite[p. 352]{LeSh93}), $(q_{\alpha,i})_{i\geq 0}$ satisfies the linear recurrence
    \begin{equation}\label{rr_qn}
    q_{\alpha,i+2s}-t_{\alpha}q_{\alpha,i+s}+(-1)^sq_{\alpha,i}=0,\ i\geq r,  
    \end{equation}
    where $t_{\alpha}$ is defined as \begin{equation}\label{t_alpha}
    t_{\alpha}=\textrm{trace} \left(\prod_{0\leq j<s}\begin{pmatrix}
    b_j & 1\\
    1 & 0
\end{pmatrix}\right).
\end{equation}\\
The next lemma gives a Binet-type formula for the convergent denominators.
\begin{lemma}\label{Binet-type}
   Let $\alpha$ be a real quadratic irrational number with simple continued fraction expansion  $[a_0;\ldots,a_{r-1},\overline{b_0,\ldots,b_{s-1}}]$, $r\geq 0, s\geq 1$. For $j=0,\ldots,s-1$, let $q^{(j)}_{\alpha,i}=q_{\alpha,si+j+r},\ i\geq 0$. 
    We have \begin{equation}\label{binet_rep}
       q^{(j)}_{\alpha,i}=c_{\alpha,1}^{(j)}\theta_{\alpha,1}^{i}-c_{\alpha,2}^{(j)}\theta_{\alpha,2}^i,\ i \geq 0,
       \end{equation}
       where 
    \begin{align*}
        \theta_{\alpha,1}&=\frac{t_{\alpha}+\sqrt{t_{\alpha}^2-4(-1)^s}}{2},\  \theta_{\alpha,2}=\frac{t_{\alpha}-\sqrt{t_{\alpha}^2-4(-1)^s}}{2},\\
        c^{(j)}_{\alpha,1}&=\frac{q^{(j)}_{\alpha,1}-\theta_{\alpha,2}q^{(j)}_{\alpha,0}}{\theta_{\alpha,1}-\theta_{\alpha,2}} > 0,\ c_{\alpha,2}^{(j)}=\frac{q^{(j)}_{\alpha,1}-\theta_{\alpha,1}q^{(j)}_{\alpha,0}}{\theta_{\alpha,1}-\theta_{\alpha,2}}.
    \end{align*}
\end{lemma}

\begin{proof}
    From \eqref{rr_qn} and the definition of $q^{(j)}_{\alpha,i}$, we have
    \[
    q^{(j)}_{\alpha,i+2}-t_{\alpha}q^{(j)}_{\alpha,i+1}+(-1)^sq^{(j)}_{\alpha,i}=0,\ i\geq 0.
    \]
Now it is easy to see that the lemma follows.
\end{proof}
\noindent Using the above lemma, we note that $1< \theta_{\alpha,1} < t_{\alpha}+1$ and $-1<\theta_{\alpha, 2} < 1$.\\

We recall that two complex numbers $u_{1}$ and $u_{2}$ are said to be \textit{multiplicatively independent} if for integers $k_{1}$ and $k_{2}$, $u_{1}^{k_{1}} u_{2}^{k_2}=1$ implies $k_{1}=k_{2}=0$.
The following lemma provides a sufficient condition for the multiplicative independence of $\theta_{\alpha,1}$ and $\theta_{\beta,1}$.
\begin{lemma}\label{TwoOst_rootsLI}
   Let $\alpha, \beta$ be real quadratic irrational numbers with simple continued fraction expansions 
    \begin{align*}
        \alpha&=[a_{0,0};\ldots,a_{0,r_0-1},\overline{b_{0,0},\ldots,b_{0,s_0-1}}], r_0\geq 0, s_0\geq 1,\\
        \beta&=[a_{1,0};\ldots,a_{1,r_1-1},\overline{b_{1,0},\ldots,b_{1,s_1-1}}], r_1\geq 0, s_1\geq 1.
    \end{align*} If $\mathbb{Q}(\alpha) \neq \mathbb{Q}(\beta)$, then $\theta_{\alpha,1}$ and $\theta_{\beta,1}$ are multiplicatively independent.
\end{lemma}
\begin{proof}
Let $(p_{\alpha, n})_{n \geq 0}$ and $(q_{\alpha,n})_{n \geq 0}$ be the sequences of convergent numerators and denominators, respectively, to the continued fraction expansion of $\alpha$.
Referring to \cite{LeSh93}, $(p_{\alpha,i})_{i\geq 0}$ satisfies the linear recurrence
    \begin{equation*}\label{rr_pn}
    p_{\alpha,i+2s}-t_{\alpha}p_{\alpha,i+s}+(-1)^s p_{\alpha,i}=0,\ i\geq r.  
    \end{equation*}    
\noindent In a similar manner to Lemma \ref{Binet-type}, we will have 
\[p^{(j)}_{\alpha,i}=d_{\alpha,1}^{(j)}\theta_{\alpha,1}^{i}-d_{\alpha,2}^{(j)}\theta_{\alpha,2}^i,\ i \geq 0,
\]
where $d^{(j)}_{\alpha,1}=\frac{p^{(j)}_{\alpha,1}-\theta_{\alpha,2}p^{(j)}_{\alpha,0}}{\theta_{\alpha,1}-\theta_{\alpha,2}} \text{ and } d_{\alpha,2}^{(j)}=\frac{p^{(j)}_{\alpha,1}-\theta_{\alpha,1}p^{(j)}_{\alpha,0}}{\theta_{\alpha,1}-\theta_{\alpha,2}}$.
We know that $ \alpha = \lim_{n \to \infty}
    \frac{p_{\alpha, n}}{q_{\alpha, n}}$. 
    Therefore, for any fixed $j$, we have
\begin{align*} 
\alpha =\lim_{i \to \infty} \frac{p_{\alpha, i}^{(j)}}{q_{\alpha, i}^{(j)}} = \lim_{i \to \infty} \frac{d_{\alpha,1}^{(j)}\theta_{\alpha,1}^{i}-d_{\alpha,2}^{(j)}\theta_{\alpha,2}^i}{c_{\alpha,1}^{(j)}\theta_{\alpha,1}^{i}-c_{\alpha,2}^{(j)}\theta_{\alpha,2}^i} =\frac{d_{\alpha,1}^{(j)}}{c_{\alpha,1}^{(j)}}
&=\frac{p^{(j)}_{\alpha,1}-\theta_{\alpha,2}p^{(j)}_{\alpha,0}}{q^{(j)}_{\alpha,1}-\theta_{\alpha,2}q^{(j)}_{\alpha,0}}.
\end{align*}
As $p^{(j)}_{\alpha,1}, p^{(j)}_{\alpha,0}, q^{(j)}_{\alpha,1},$ and $q^{(j)}_{\alpha,0}$ are all integers, we obtain that $\mathbb{Q}(\alpha) = \mathbb{Q}(\theta_{\alpha, 2}) = \mathbb{Q}(\sqrt{t_{\alpha}^{2}-4(-1)^{s_0}})= \mathbb{Q}(\theta_{\alpha, 1}).$
As $\mathbb{Q}(\alpha) \neq \mathbb{Q}(\beta)$, we have $\mathbb{Q}(\theta_{\alpha,1}) \neq \mathbb{Q}(\theta_{\beta,1})$. Therefore, $\theta_{\alpha,1}$ and $\theta_{\beta,1}$ have to be multiplicatively independent.
\end{proof}

\begin{lemma}\label{diff_ublb}
    Let $\alpha$ be a real quadratic irrational number with simple continued fraction expansion $[a_0;\ldots,a_{r-1},\overline{b_0,\ldots,b_{s-1}}]$, $r\geq 0, s\geq 1$. Then 
     \begin{equation*}
         C_{\alpha, 3} \theta_{\alpha,1}^{n_1} \leq \left| q_{\alpha, N_1} - q_{\alpha, N_2}\right| \leq C_{\alpha, 4} \theta_{\alpha,1}^{n_1}.
     \end{equation*}
\end{lemma}
\begin{proof}
  Recall that $q^{(j)}_{\alpha,i}=q_{\alpha,si+j+r}$, from this we obtain $q_{\alpha, N}=q_{\alpha, n}^{(j)}$.
Using equation \eqref{binet_rep}, we have
\[
q_{\alpha, n}^{(j)} = c_{\alpha, 1}^{(j)}\theta_{\alpha, 1}^{n} - c_{\alpha, 2}^{(j)}\theta_{\alpha, 2}^{n}. 
\]
As $n \to \infty$ \[ q_{\alpha, n}^{(j)}  \sim c_{\alpha, 1}^{(j)}\theta_{\alpha, 1}^{n}.\] 
Hence,  \[  q_{\alpha, N} =  q_{\alpha, n}^{(j)} \sim c_{\alpha, 1}^{(j)} \theta_{\alpha, 1}^{n},\ \ \text{ as } N \to \infty.\]
Therefore, there are positive constants $C_{\alpha, 5}$ and $C_{\alpha, 6}$ such that $\nicefrac{C_{\alpha, 6}}{C_{\alpha, 5}} < \theta_{\alpha,1}$ with
\[ C_{\alpha, 5} \theta_{\alpha,1}^{n} \leq  q_{\alpha, N} \leq C_{\alpha, 6} \theta_{\alpha,1}^{n}, \ \ N\geq N_0,\] where $N_{0}$ is a sufficiently large integer.\\
For $N_1 > N_2 \geq N_0$, 
\begin{align}
\nonumber
  \left| q_{\alpha, N_1} - q_{\alpha, N_2}\right| \leq \left|q_{\alpha, N_1} \right| + \left|q_{\alpha, N_2} \right| &\leq C_{\alpha, 6} (\theta_{\alpha,1}^{n_1}+\theta_{\alpha,1}^{n_2})\\
\nonumber  &= C_{\alpha, 6} \left(1+\frac{1}{\theta_{\alpha,1}^{n_1 - n_2}}\right) \theta_{\alpha,1}^{n_1}\\
\label{ub_1}
    &\leq C_{\alpha, 4} \theta_{\alpha,1}^{n_1}, 
\end{align}
 where  $C_{\alpha, 4} =2C_{\alpha, 6}.$\\
For $N_1 > N_2 \geq N_0$, with $n_1 \neq n_2$, we have
\begin{align}
\nonumber
    \left| q_{\alpha, N_1} - q_{\alpha, N_2}\right| \geq \left|q_{\alpha, N_1} \right| - \left|q_{\alpha, N_2} \right|
    &\geq C_{\alpha, 5}\theta_{\alpha,1}^{n_1}- C_{\alpha, 6} \theta_{\alpha,1}^{n_2}\\
\nonumber
    &= C_{\alpha, 5}  \left(1-\frac{C_{\alpha, 6}}{C_{\alpha, 5} \theta_{\alpha,1}^{n_1 - n_2}}\right) \theta_{\alpha,1}^{n_1}\\
\label{lb_1}
    &\geq C_{\alpha, 5} \left( 1-\frac{C_{\alpha, 6}}{C_{\alpha, 5} \theta_{\alpha,1}}\right) \theta_{\alpha,1}^{n_1},
\end{align}
and with $n_1 = n_2$, using equation \eqref{binet_rep}, we have

\begin{align*}
    \left| q_{\alpha, N_1} - q_{\alpha, N_2}\right| &= \left| (c_{\alpha,1}^{(j_1)} - c_{\alpha,1}^{(j_2)}) \theta_{\alpha,1}^{n_1} - (c_{\alpha,2}^{(j_1)} - c_{\alpha,2}^{(j_2)}) \theta_{\alpha,2}^{n_1} \right|\\
    &= \left| (c_{\alpha,1}^{(j_1)} - c_{\alpha,1}^{(j_2)})  - (c_{\alpha,2}^{(j_1)} - c_{\alpha,2}^{(j_2)}) \left(\frac{\theta_{\alpha,2}}{\theta_{\alpha,1}}\right)^{n_1} \right| \theta_{\alpha,1}^{n_1}.    
\end{align*}

\noindent Since $\lim_{n_1 \to \infty}  \left(\frac{\theta_{\alpha,2}}{\theta_{\alpha,1}}\right)^{n_1} = 0$, there exists $n^{*} \in \mathbb{N}$ such that
\[ \left| \left(\frac{\theta_{\alpha,2}}{\theta_{\alpha,1}}\right)^{n_1} \right| \leq  \frac{|c_{\alpha,1}^{(j_1)} - c_{\alpha,1}^{(j_2)}|}{2 |c_{\alpha,2}^{(j_1)} - c_{\alpha,2}^{(j_2)}|},\ \ \ n_1 \geq n^{*}, \]
therefore
\begin{align}
\label{lb_2}
    \left| q_{\alpha, N_1} - q_{\alpha, N_2}\right| 
    & \geq  \frac{\left|c_{\alpha,1}^{(j_1)} - c_{\alpha,1}^{(j_2)}  \right|}{2} \theta_{\alpha,1}^{n_1}, 
\end{align}
for $ N_1 > N_2 \geq \max (N_0, N_{0}^{*} = s n^{*} + j_2 + r).$\\
\noindent Combining equations \eqref{lb_1}
and \eqref{lb_2}, we have
\begin{equation}\label{lb_3}
\left| q_{\alpha, N_1} - q_{\alpha, N_2}\right| \geq C_{\alpha, 3} \theta_{\alpha,1}^{n_1},    
\end{equation}
 where $C_{\alpha, 3} = \min\left(C_{\alpha, 5} \left( 1-\frac{C_{\alpha, 6}}{C_{\alpha, 5} \theta_{\alpha,1}}\right), \min_{0 \leq j_2 < j_1 \leq s-1} \frac{\left| c_{\alpha,1}^{(j_1)} - c_{\alpha,1}^{(j_2)}  \right|} {2} \right)$.

\noindent Combining inequalities \eqref{ub_1} and \eqref{lb_1}, we get the desired result.
\end{proof}
Next, we recall the definition and some basic properties of the absolute logarithmic height function $h(\cdot)$. If $\delta$ is an algebraic  number with minimal polynomial $f(X)=d_0(X-\delta^{(1)})\cdots(X-\delta^{(d)})\in\mathbb{Z}[X]$, 
    \[
    h(\delta)=\frac{1}{d}\left(\log d_0+\sum_{i=1}^d\max(0,\log|\delta^{(i)}|)\right)
    \]
denotes its absolute logarithmic height.
If $\delta_1,\delta_2$ are algebraic numbers, then
    \begin{align}\label{height_of_sum}
       h(\delta_1 \pm \delta_2)&\leq h(\delta_1)+h(\delta_2)+\log 2,\\ \label{height_of_product}
       h(\delta_1\delta_2^{\pm 1})&\leq h(\delta_1)+h(\delta_2),\\
       \label{height_of_exponents}
       h(\delta_1^m)&=|m|h(\delta_1)\ (m\in\mathbb{Z})
    \end{align}
(see \cite[Property 3.3]{Wa_00}). 


The following lemma due to Chim, Pink \& Ziegler \cite{Chim_Ziegler_18} gives a lower bound for the absolute logarithmic height of the product of powers of multiplicatively independent algebraic numbers.
\begin{lemma}\label{CPZ}
    Let $K$ be a number field and suppose that $\delta_1, \delta_2 \in K$ are multiplicatively independent. There exists an effectively computable constant $D_{0} > 0$ such that for $n, m \in \mathbb{Z}$, we have
    $$ h\left(\frac{\delta_{1}^{n}}{\delta_{2}^{m}}\right) \geq D_{0} \max (|n|, |m|).$$
\end{lemma}

We will use the following lower bound for linear forms in logarithms due to Matveev \cite{Mat1, Mat2}. (See also \cite{Baker_Wustholz_07, Bugead_Mignotte_Siksek_06}.)  
\begin{lemma}\label{lfl}
    Let $T$ be a positive integer. Let  $\delta_1,\ldots,\delta_T$ be non-zero real algebraic numbers and $\log\delta_1,\ldots,\log\delta_T$ be some determinations of their complex logarithms. Let $D$ be the degree of the number field generated by $\delta_1,\ldots,\delta_T$ over $\mathbb{Q}$.  For $j=1,\ldots,T$, let $A_j'$ be a real number satisfying
    \[
    A_j'\geq\max\left(Dh(\delta_j),|\log\delta_j|,0.16\right).
    \]
    Let $k_1,\ldots,k_T$ be rational integers. Set
    \[
    B=\max(|k_1|,\ldots,|k_T|) \textrm{ and }\ \Lambda=\delta_1^{k_1}\cdots\delta_T^{k_T}-1.
    \]
    If $\Lambda \neq 0$, then
    \[
    \log|\Lambda|>-1.4\times 30^{T+3}T^{4.5}D^{2}\log(eD) A_1'\cdots A_T'\log(eB).
    \]
\end{lemma}
We shall need the following lemma due to Peth\H{o} \& de Weger \cite{Petho_Weger_86}.
\begin{lemma}\label{PW}
    Let $a \geq 0$, $c \geq 1$, $g> (\nicefrac{e^2}{c})^{c}$, and let $x \in \mathbb{R}$ be the largest solution of $x= a+g(\log x)^{c}.$ Then, 
    $$ x < 2^{c} (a^{1/c} + g^{1/c} \log (c^{c}g))^{c}.$$
\end{lemma}

The next result is a slight variation of a result due to Dujella and Peth\H{o} \cite{DP_98}, which itself is a generalization of the result due to Baker and Davenport \cite{BD_69}. For $x \in \mathbb{R}$, let $\left\|x\right\| = \min \left( x-  \left \lfloor x \right \rfloor , \left \lceil x \right \rceil -x \right)$ be the distance from $x$ to the nearest integer.
\begin{lemma}\label{BD_reduction}
    Let $M$ be a positive integer. Let $q$ be a convergent denominator to the continued fraction of the irrational $\gamma$ such that $q > 6M$, and let $A_1, A_2, \kappa$ be real numbers with $A_1>0$ and $A_2 >1$. Let $\varepsilon:= \left\| \kappa q\right\| - M\left\|\gamma q\right\|.$ If $\varepsilon >0$, then there is no solution to the inequality \[ 0< m \gamma - n+\kappa < A_1 A_2^{-k}\] in positive integers $m, n$ and $k$ with $m \leq M$ and $k \geq \dfrac{\log (\nicefrac{A_{1} q}{\varepsilon})}{\log A_2}$.
\end{lemma}
\section{Proof of Theorem \ref{main thm}}\label{Sec_main thm}

In what follows, $C_1,C_2,\ldots$ denote constants depending only upon\\ $s_0,s_1,a_{0,0},\ldots,a_{0,r_0-1},b_{0,0},\ldots,b_{0,s_0-1},a_{1,0},\ldots,a_{1,r_1-1},b_{1,0},\ldots,b_{1,s_1-1}$.

 Let $c$ be an integer. Assume that $(N_1, M_1) \neq (N_2, M_2)$ are pairs of indices with $N_{2} \leq N_{1}$ such that 
\begin{equation}\label{dist_rep_c}
c = q_{\alpha, N_{1}}-q_{\beta, M_{1}} = q_{\alpha, N_{2}}- q_{\beta, M_{2}},
\end{equation}
i.e.
\begin{equation}\label{main_eqn}
 q_{\alpha, N_{1}}- q_{\alpha, N_{2}} =q_{\beta, M_{1}} - q_{\beta, M_{2}}.
\end{equation}
Since the sequences $\left(q_{\alpha, i}\right)_{i \geq 0}$, $\left(q_{\beta, i}\right)_{i \geq 0}$ are either strictly increasing or would eventually be strictly increasing, there exists $N^{*}, M^{*} \in \mathbb{N}$ such that for $i_{1} > i_{2} \geq N^{*}$ and $j_{1} > j_{2} \geq M^{*},$ $q_{\alpha, i_1} > q_{\alpha, i_2}$ and $q_{\beta, j_1} > q_{\beta, j_2},$ respectively.
For $N_{1} < N^{*}$, $c$ has finitely many possibilities. So, from now on, $N_{1} \geq N^{*}$.\\ 
\indent We further assume that $N_1 \neq N_2$, since otherwise we either have $M_{1} = M_{2}$ or $M_1 \neq M_2$ with both $M_1, M_2 \leq M^{*}$ which would give finite choices for $c$. Without loss of generality, we assume $N_1 > N_2$, which implies $q_{\alpha, N_{1}}- q_{\alpha, N_{2}} > 0$. Using equation \eqref{main_eqn}, $q_{\beta, M_{1}} - q_{\beta, M_{2}} >0$ and therefore $M_1 > M_2$.\\
 Write $N_i=s_0n_i+j_i+r_0$, $0\leq j_i\leq s_0-1$ and
$M_i=s_1m_i+p_i+r_1$, $0\leq p_i\leq s_1-1$ for $i=1,2$. For technical reasons, we assume that $N_1 \geq 2s_{0}+r_{0}$, so that $n_1 \geq 2$. In fact, for $N_1 <  2 s_{0} + r_{0}$, there will only be finitely many possibilities for $c$. \\

\noindent Applying Lemma \ref{diff_ublb} for $\beta$, we get

\begin{equation}\label{ub_lb}
    C_{\beta, 3} \theta_{\beta,1}^{m_1} \leq \left| q_{\beta, M_1} - q_{\beta, M_2}\right| \leq C_{\beta, 4} \theta_{\beta,1}^{m_1}. 
\end{equation} 
Using equations \eqref{ub_1}, \eqref{lb_3}, \eqref{main_eqn}, and \eqref{ub_lb}, we get

\begin{equation}\label{domterm_samesize}
   C_{\alpha, 3} \theta_{\alpha,1}^{n_1} \leq \left| q_{\alpha, N_1} - q_{\alpha, N_2}\right|= \left| q_{\beta, M_1} - q_{\beta, M_2}\right| \leq C_{\beta, 4} \theta_{\beta,1}^{m_1}
\end{equation}
and
\begin{equation}\label{domterm_samesize2}
   C_{\beta, 3} \theta_{\beta,1}^{m_1} \leq \left| q_{\beta, M_1} - q_{\beta, M_2}\right|= \left| q_{\alpha, N_1} - q_{\alpha, N_2}\right| \leq C_{\alpha, 4} \theta_{\alpha,1}^{n_1}. 
\end{equation}
Using \eqref{domterm_samesize2}, we get
\begin{equation}\label{LFL1_B}
        m_1<C_{1} n_1,\ \textrm{ where } \ C_{1} = \dfrac{\log \theta_{\alpha, 1}+2 \log^+ (\nicefrac{C_{\alpha, 4}}{C_{\beta, 3}})}{\log \theta_{\beta, 1}}
    \end{equation}
since $n_1 \geq 2$ and $\log^+(x)=\max(\log x,0)$.\\
\noindent Now, using equations \eqref{binet_rep} and \eqref{main_eqn}, we get
\begin{align}
\nonumber
 c_{\alpha, 1}^{(j_1)}\theta_{\alpha, 1}^{n_1} - c_{\alpha, 2}^{(j_1)}\theta_{\alpha, 2}^{n_1} - c_{\alpha, 1}^{(j_2)}\theta_{\alpha, 1}^{n_2} + c_{\alpha, 2}^{(j_2)}\theta_{\alpha, 2}^{n_2}
&= 
c_{\beta, 1}^{(p_1)} \theta_{\beta, 1}^{m_1} - c_{\beta, 2}^{(p_1)} \theta_{\beta, 2}^{m_1} - c_{\beta, 1}^{(p_2)} \theta_{\beta, 1}^{m_2} \\ \label{main_eqn_1} &+ c_{\beta, 2}^{(p_2)} \theta_{\beta, 2}^{m_2}.   
\end{align}

\noindent There are $4$ possible cases:
\begin{enumerate}[(i)]
    \item $n_{1} = n_{2}$ and $m_{1} = m_{2}$;
    \item $n_{1} = n_{2}$ and $m_{1} > m_{2}$;
    \item $n_{1} > n_{2}$ and $m_{1} = m_{2}$;
    \item $n_{1} > n_{2}$ and $m_{1} > m_{2}$.
\end{enumerate}
In each of these cases, we show that $m_{1} \ll (\log n_1)^{3}$, where $\ll$ indicates that the constant involved depends only on $\alpha$ and $\beta$.\\
\noindent \textbf{Case (i):} $n_{1} = n_{2}$ and $m_{1} = m_{2}$.\\
Using equation \eqref{main_eqn_1}, we have
 \[(c_{\alpha, 1}^{(j_1)} - c_{\alpha, 1}^{(j_2)}) \theta_{\alpha, 1}^{n_1} - (c_{\alpha, 2}^{(j_1)} - c_{\alpha, 2}^{(j_2)}) \theta_{\alpha, 2}^{n_1}  = 
(c_{\beta, 1}^{(p_1)} - c_{\beta, 1}^{(p_2)}) \theta_{\beta, 1}^{m_1} - (c_{\beta, 2}^{(p_1)} - c_{\beta, 2}^{(p_2)}) \theta_{\beta, 2}^{m_1},\]
i.e.
 \[(c_{\alpha, 1}^{(j_1)} - c_{\alpha, 1}^{(j_2)}) \theta_{\alpha, 1}^{n_1} - (c_{\beta, 1}^{(p_1)} - c_{\beta, 1}^{(p_2)}) \theta_{\beta, 1}^{m_1} = (c_{\alpha, 2}^{(j_1)} - c_{\alpha, 2}^{(j_2)}) \theta_{\alpha, 2}^{n_1}- (c_{\beta, 2}^{(p_1)} - c_{\beta, 2}^{(p_2)}) \theta_{\beta, 2}^{m_1}.\]
 
\noindent Since $|\theta_{\alpha, 2}|<1$, $|\theta_{\beta, 2}|<1$, we get\\
\noindent
 \[\left|(c_{\alpha, 1}^{(j_1)} - c_{\alpha, 1}^{(j_2)})  \theta_{\alpha, 1}^{n_1} - (c_{\beta, 1}^{(p_1)} - c_{\beta, 1}^{(p_2)})  \theta_{\beta, 1}^{m_1}  \right| \leq C_{2}.\]

\noindent  If $c_{\beta, 1}^{(p_1)} - c_{\beta, 1}^{(p_2)} = 0$, then, using Lemma \ref{Binet-type}, we obtain
\[ q_{\beta, 1}^{(p_1)} - q_{\beta, 1}^{(p_2)} = \theta_{\beta, 2} (q_{\beta, 0}^{(p_1)} - q_{\beta, 0}^{(p_2)}).\] We observe that the left-hand side is an integer, while the right-hand side is an irrational. Therefore, $p_1 = p_2$. We recall that $M_i = s_1 m_i + p_i + r_1$, this gives $M_1 = M_2$, which contradicts our assumption that $M_1 > M_2$. Hence, $c_{\beta, 1}^{(p_1)} - c_{\beta, 1}^{(p_2)} \neq 0$.\\ 
 Dividing by $\left|(c_{\beta, 1}^{(p_1)} - c_{\beta, 1}^{(p_2)})  \theta_{\beta, 1}^{m_1}\right|$, we get 

 \begin{align}
     \nonumber \left|\dfrac{c_{\alpha, 1}^{(j_1)} - c_{\alpha, 1}^{(j_2)}}{c_{\beta, 1}^{(p_1)} - c_{\beta, 1}^{(p_2)}}  \theta_{\alpha, 1}^{n_1} \theta_{\beta, 1}^{-m_1}- 1    \right| &\leq \frac{C_{2}}{\left|c_{\beta, 1}^{(p_1)} - c_{\beta, 1}^{(p_2)}\right|  \theta_{\beta, 1}^{m_1}}\\
\label{lfl1_rhs}
     &=\frac{C_{3}}{\theta_{\beta, 1}^{m_1}}, 
 \end{align}
 where $C_{3} = \dfrac{C_{2}}{\left|c_{\beta, 1}^{(p_1)} - c_{\beta, 1}^{(p_2)}\right|}$.
 
\noindent Let 
    \[
    \Gamma_{1}=\frac{c_{\alpha, 1}^{(j_1)} - c_{\alpha, 1}^{(j_2)}}{c_{\beta, 1}^{(p_1)} - c_{\beta, 1}^{(p_2)}}  \theta_{\alpha, 1}^{n_1} \theta_{\beta, 1}^{-m_1}.
    \]
If $\Gamma_{1} =1$, then $\dfrac{c_{\alpha, 1}^{(j_1)} - c_{\alpha, 1}^{(j_2)}}{c_{\beta, 1}^{(p_1)} - c_{\beta, 1}^{(p_2)}}=\dfrac{\theta_{\beta, 1}^{m_1}}{\theta_{\alpha, 1}^{n_1}}$. Taking heights and applying Lemma \ref{CPZ}, we get
\begin{equation}\label{case(i)_bounded m_1 1}
    \max (n_{1}, m_{1}) \ll 1.
\end{equation}

\noindent If $\Gamma_{1}\neq 1$, then we apply Lemma \ref{lfl} with $T=3$, $D=4$,
    \begin{align*}
        \delta_1&=\theta_{\alpha, 1},\ \delta_2=\theta_{\beta, 1},\ \delta_3=\dfrac{c_{\alpha, 1}^{(j_1)} - c_{\alpha, 1}^{(j_2)}}{c_{\beta, 1}^{(p_1)} - c_{\beta, 1}^{(p_2)}},\\
        k_1&=n_1,\ k_2=-m_1,\ k_3=1,\\
        A_1'&\ll 1,\ A_2'\ll 1,\ A_3' \ll 1.
    \end{align*}
Using \eqref{LFL1_B}, we can take $B=\max(C_{1},1)n_1$. Then 
    \[
  \log|\Gamma_{1}-1|>-C_{4}\log(e \max(C_{1},1)n_1).  
    \]
Comparing with \eqref{lfl1_rhs}, we obtain
\[ m_{1}< \frac{C_{4}(1+\log \max(C_1, 1)+ \log n_{1})+\log C_{3}}{\log \theta_{\beta,1}}.\] Therefore, 
    \begin{equation}\label{case(i)_bounded m_1 2} 
    m_{1}<C_{5}\log n_1,
    \end{equation}
where $C_{5}=\dfrac{2 C_{4}\left(2+\log^+ C_1 \right)+2\log^{+} C_{3}}{\log \theta_{\beta,1}}$ since $n_1 \geq 2$.

\noindent From inequalities \eqref{case(i)_bounded m_1 1} and \eqref{case(i)_bounded m_1 2}, we obtain
$$ m_1 \ll \log n_1 \ll (\log n_1)^{3}.$$

\noindent \textbf{Case (ii):} $n_{1} = n_{2}$ and $m_{1} > m_{2}$.\\
Using equation \eqref{main_eqn_1}, we have
\begin{equation}\label{scII main_eqn}
    (c_{\alpha, 1}^{(j_1)} - c_{\alpha, 1}^{(j_2)}) \theta_{\alpha, 1}^{n_1} - (c_{\alpha, 2}^{(j_1)} - c_{\alpha, 2}^{(j_2)}) \theta_{\alpha, 2}^{n_1}  = 
c_{\beta, 1}^{(p_1)} \theta_{\beta, 1}^{m_1} - c_{\beta, 2}^{(p_1)} \theta_{\beta, 2}^{m_1} - c_{\beta, 1}^{(p_2)} \theta_{\beta, 1}^{m_2} + c_{\beta, 2}^{(p_2)} \theta_{\beta, 2}^{m_2},
\end{equation}
i.e.
 \[(c_{\alpha, 1}^{(j_1)} - c_{\alpha, 1}^{(j_2)}) \theta_{\alpha, 1}^{n_1} - c_{\beta, 1}^{(p_1)} \theta_{\beta, 1}^{m_1} =  (c_{\alpha, 2}^{(j_1)} - c_{\alpha, 2}^{(j_2)}) \theta_{\alpha, 2}^{n_1} - c_{\beta, 2}^{(p_1)} \theta_{\beta, 2}^{m_1} - c_{\beta, 1}^{(p_2)} \theta_{\beta, 1}^{m_2} + c_{\beta, 2}^{(p_2)} \theta_{\beta, 2}^{m_2}.\]
Since $|\theta_{\alpha, 2}|<1$, $|\theta_{\beta, 2}|<1$, we get
\[ \left|(c_{\alpha, 1}^{(j_1)} - c_{\alpha, 1}^{(j_2)}) \theta_{\alpha, 1}^{n_1} - c_{\beta, 1}^{(p_1)} \theta_{\beta, 1}^{m_1} \right| \leq \left|c_{\alpha, 2}^{(j_1)} - c_{\alpha, 2}^{(j_2)}\right| +  \left|c_{\beta, 2}^{(p_1)}\right|  +  c_{\beta, 1}^{(p_2)}\theta_{\beta, 1}^{m_2}+ \left|c_{\beta, 2}^{(p_2)}\right|.\]
Dividing by $c_{\beta, 1}^{(p_1)}  \theta_{\beta, 1}^{m_1}$, we get
\begin{align}
\nonumber
   \left| \frac{c_{\alpha, 1}^{(j_1)} - c_{\alpha, 1}^{(j_2)}}{c_{\beta, 1}^{(p_1)}} \theta_{\alpha, 1}^{n_1} \theta_{\beta, 1}^{-m_1}- 1  \right| &\leq \frac{\left|c_{\alpha, 2}^{(j_1)} - c_{\alpha, 2}^{(j_2)}\right| +  \left|c_{\beta, 2}^{(p_1)}\right|  +   \left|c_{\beta, 2}^{(p_2)}\right|}{c_{\beta, 1}^{(p_1)}  \theta_{\beta, 1}^{m_1}} +\frac{c_{\beta, 1}^{(p_2)}\theta_{\beta, 1}^{m_2}}{c_{\beta, 1}^{(p_1)}  \theta_{\beta, 1}^{m_1}}\\
\nonumber
 &\leq \frac{C_{6}}{\theta_{\beta, 1}^{m_1}} + \frac{C_{7}}{\theta_{\beta, 1}^{m_1-m_2}}\\
 \label{lfl2_rhs}
 &\leq \frac{C_{8}}{\theta_{\beta, 1}^{m_1-m_2}}.
\end{align}

\noindent Let 
    \[
    \Gamma_{2}=\frac{c_{\alpha, 1}^{(j_1)} - c_{\alpha, 1}^{(j_2)}}{c_{\beta, 1}^{(p_1)}} \theta_{\alpha, 1}^{n_1} \theta_{\beta, 1}^{-m_1}.
    \]
If $\Gamma_{2} =1$, then $\dfrac{c_{\alpha, 1}^{(j_1)} - c_{\alpha, 1}^{(j_2)}}{c_{\beta, 1}^{(p_1)}}=\dfrac{\theta_{\beta, 1}^{m_1}}{\theta_{\alpha, 1}^{n_1}}$. Taking heights and applying Lemma \ref{CPZ}, we get
\begin{equation}\label{case(ii)_bounded m_1 1}
  \max (n_{1}, m_{1}) \ll 1.  
\end{equation}

\noindent If $\Gamma_{2}\neq 1$, then we apply Lemma \ref{lfl} with $T=3$, $D=4$,
    \begin{align*}
        \delta_1&=\theta_{\alpha, 1},\ \delta_2=\theta_{\beta, 1},\ \delta_3=\dfrac{c_{\alpha, 1}^{(j_1)} - c_{\alpha, 1}^{(j_2)}}{c_{\beta, 1}^{(p_1)}},\\
        k_1&=n_1,\ k_2=-m_1,\ k_3=1,\\
        A_1'&\ll 1,\ A_2'\ll 1,\ A_3'\ll 1.
    \end{align*}
Using \eqref{LFL1_B}, we can take $B=\max(C_{1},1)n_1$. Then 
    \[
  \log|\Gamma_{2}-1|>-C_{9}\log(e \max(C_{1},1)n_1).  
    \]
Comparing with \eqref{lfl2_rhs}, we obtain
\[ m_{1} - m_{2}< \frac{C_{9}(1+\log \max(C_1, 1)+ \log n_{1})+\log C_{8}}{\log \theta_{\beta,1}}.\] Therefore,
    \begin{equation*}
    m_{1} -m_{2}< C_{10} \log n_1,
    \end{equation*}
where $C_{10}=\dfrac{2 C_{9}\left(2+\log^+ C_1 \right)+2 \log^{+} C_{8}}{\log \theta_{\beta,1}}$ since $n_1 \geq 2$.\\

Again using equation \eqref{scII main_eqn}, we have
\[(c_{\alpha, 1}^{(j_1)} - c_{\alpha, 1}^{(j_2)}) \theta_{\alpha, 1}^{n_1} - c_{\beta, 1}^{(p_1)} \theta_{\beta, 1}^{m_1} + c_{\beta, 1}^{(p_2)} \theta_{\beta, 1}^{m_2} =  (c_{\alpha, 2}^{(j_1)} - c_{\alpha, 2}^{(j_2)}) \theta_{\alpha, 2}^{n_1} - c_{\beta, 2}^{(p_1)} \theta_{\beta, 2}^{m_1} + c_{\beta, 2}^{(p_2)} \theta_{\beta, 2}^{m_2}.\]

\noindent Since $|\theta_{\alpha, 2}|<1$, $|\theta_{\beta, 2}|<1$, we get

\[ \left|(c_{\alpha, 1}^{(j_1)} - c_{\alpha, 1}^{(j_2)}) \theta_{\alpha, 1}^{n_1} - c_{\beta, 1}^{(p_1)} \theta_{\beta, 1}^{m_1} + c_{\beta, 1}^{(p_2)} \theta_{\beta, 1}^{m_2} \right| \leq \left|(c_{\alpha, 2}^{(j_1)} - c_{\alpha, 2}^{(j_2)}) \right| + \left|c_{\beta, 2}^{(p_1)} \right|   + \left|c_{\beta, 2}^{(p_2)} \right|. \]

\noindent Dividing by $\left|c_{\beta, 1}^{(p_1)}\theta_{\beta, 1}^{m_1} - c_{\beta, 1}^{(p_2)}\theta_{\beta, 1}^{m_2}\right|$, we get

\begin{align*}
    \left| \frac{c_{\alpha, 1}^{(j_1)} - c_{\alpha, 1}^{(j_2)}}{c_{\beta, 1}^{(p_1)} - c_{\beta, 1}^{(p_2)}\theta_{\beta, 1}^{m_2-m_1}} \theta_{\alpha, 1}^{n_1} \theta_{\beta, 1}^{-m_1} - 1 \right| 
    &\leq  \frac{\left|(c_{\alpha, 2}^{(j_1)} - c_{\alpha, 2}^{(j_2)}) \right| + \left|c_{\beta, 2}^{(p_1)} \right|   + \left|c_{\beta, 2}^{(p_2)} \right|} {\left|c_{\beta, 1}^{(p_1)} - c_{\beta, 1}^{(p_2)}\theta_{\beta, 1}^{-(m_1 - m_2)} \right| } \frac{1}{\theta_{\beta, 1}^{m_1}}.
\end{align*}

We now find an upper bound for the expression $\dfrac{1}{\left|c_{\beta, 1}^{(p_1)} - c_{\beta, 1}^{(p_2)}\theta_{\beta, 1}^{-(m_1 - m_2)} \right|}.$
The sequence $(\frac{1}{\theta_{\beta,1}^{k}})_{k \geq 1}$ decreases to $0$. Therefore, there exists a $k_{0}\in \mathbb{N}$ depending only on $\beta$ such that \[ \left| \frac{1}{\theta_{\beta,1}^{k}}  \right| < \frac{c_{\beta, 1}^{(p_1)}}{2c_{\beta, 1}^{(p_2)}} , \text { for all } k \geq k_0.\] 

\noindent Since $c_{\beta, 1}^{(p_1)},  c_{\beta, 1}^{(p_2)}$ and $\theta_{\beta,1}$ are positive numbers, we have \[
  \left|c_{\beta, 1}^{(p_1)} - c_{\beta, 1}^{(p_2)}\theta_{\beta, 1}^{-k} \right| > c_{\beta, 1}^{(p_1)}  -  c_{\beta, 1}^{(p_2)}\theta_{\beta, 1}^{-k} > \frac{c_{\beta, 1}^{(p_1)}}{2}.
\]

\noindent Now we will show that for $1 \leq k \leq k_{0}-1$, $\left|c_{\beta, 1}^{(p_1)} - c_{\beta, 1}^{(p_2)}\theta_{\beta, 1}^{-k} \right| \neq 0$.\\
\noindent Suppose $\left|c_{\beta, 1}^{(p_1)} - c_{\beta, 1}^{(p_2)}\theta_{\beta, 1}^{-k} \right| = 0$. Substituting the values of $c_{\beta, 1}^{(p_1)}, c_{\beta, 1}^{(p_2)}$ from Lemma \ref{Binet-type}, we get
\[ \theta_{\beta, 1}^{k} q_{\beta, 1}^{(p_1)}- \theta_{\beta, 1}^{k} \theta_{\beta, 2} q_{\beta, 0}^{(p_1)} + \theta_{\beta, 2} q_{\beta, 0}^{(p_2)} = q_{\beta, 1}^{(p_2)}.\]
Since $\theta_{\beta, 1} \theta_{\beta, 2} = (-1)^{s_1}$, we have

\begin{equation}\label{case(ii)_non-zero exp}
    \theta_{\beta, 1}^{k} q_{\beta, 1}^{(p_1)}- (-1)^{s_1} \theta_{\beta, 1}^{k-1}  q_{\beta, 0}^{(p_1)} + \theta_{\beta, 2} q_{\beta, 0}^{(p_2)} = q_{\beta, 1}^{(p_1)}.
\end{equation} 

\noindent Let $\theta_{\beta, 1}^{k-1} = A_{k-1} \theta_{\beta, 1} + B_{k-1}$, where $A_{k-1}, B_{k-1}$ are integers. Multiplying both sides by $\theta_{\beta, 1}$ and using the minimal polynomial of $\theta_{\beta, 1}$, we get $\theta_{\beta, 1}^{k} = (A_{k-1} t_{\beta} + B_{k-1})\theta_{\beta, 1} -(-1)^{s_1} A_{k-1}.$
Substituting the expressions for $\theta_{\beta, 1}^{k}, \theta_{\beta, 1}^{k-1}$ into equation \eqref{case(ii)_non-zero exp} and using equation \eqref{rr_qn}, we get

\begin{align*}
    A_{k-1} \theta_{\beta, 1} q_{\beta, 2}^{(p_1)} + B_{k-1} \theta_{\beta, 1} q_{\beta, 1}^{(p_1)} + \theta_{\beta, 2} q_{\beta, 0}^{(p_2)} &= q_{\beta, 1}^{(p_1)}+ (-1)^{s_1}A_{k-1}q_{\beta, 1}^{(p_1)}\\ &+ (-1)^{s_1}B_{k-1}q_{\beta, 0}^{(p_1)}
\end{align*} 

 \noindent Observe that the right-hand side is an integer. We consider the irrational part that arises from the left-hand side and equate it to $0$.\\
 This gives, 
 \[ A_{k-1} q_{\beta, 2}^{(p_1)}+ B_{k-1}q_{\beta, 1}^{(p_1)}- q_{\beta, 0}^{(p_2)} = 0.\]
 
\noindent For $k=1$, we have $A_{0}=0, B_{0}=1$, this gives $q_{\beta, 1}^{(p_1)}- q_{\beta, 0}^{(p_2)} = 0,$ which is not possible as the LHS is always positive.\\  
\noindent For $k \geq 2$, we shall use a different argument.\\
 When $s_1$ is odd, we have $A_{k-1}$ and $B_{k-1}$ are non-negative integers (simultaneously not zero) and this gives positive LHS, which is a contradiction.\\
 \noindent When $s_1$ is even, it can be shown by induction that $A_{k-1} \geq |B_{k-1}|+1$. As, $B_{k-1} \leq 0$, we have
 \[ (A_{k-1}-|B_{k-1}|-1) q_{\beta, 2}^{(p_1)}+ |B_{k-1}|(q_{\beta, 2}^{(p_1)} - q_{\beta, 1}^{(p_1)})+ (q_{\beta, 2}^{(p_1)} - q_{\beta, 0}^{(p_2)}) = 0,\] which is not possible as the LHS is positive.\\
 
 \noindent Hence, for $1 \leq k \leq k_{0}-1$, $\left|c_{\beta, 1}^{(p_1)} - c_{\beta, 1}^{(p_2)}\theta_{\beta, 1}^{-k} \right| \neq 0.$\\
\noindent Let $ m := \min\left(\left|c_{\beta, 1}^{(p_1)} - c_{\beta, 1}^{(p_2)}\theta_{\beta, 1}^{-1} \right|, \ldots, \left|c_{\beta, 1}^{(p_1)} - c_{\beta, 1}^{(p_2)}\theta_{\beta, 1}^{-(k_0 - 1)}\right|, \dfrac{c_{\beta, 1}^{(p_1)}}{2}\right).$\\ 
\noindent Finally, we get
\begin{equation}\label{lfl3_rhs}
    \left| \frac{c_{\alpha, 1}^{(j_1)} - c_{\alpha, 1}^{(j_2)}}{c_{\beta, 1}^{(p_1)} - c_{\beta, 1}^{(p_2)}\theta_{\beta, 1}^{m_2-m_1}} \theta_{\alpha, 1}^{n_1} \theta_{\beta, 1}^{-m_1} - 1 \right| \leq \frac{C_{11}}{\theta_{\beta, 1}^{m_1}}.
\end{equation}

\noindent Let 
    \[
    \Gamma_{3}= \frac{c_{\alpha, 1}^{(j_1)} - c_{\alpha, 1}^{(j_2)}}{c_{\beta, 1}^{(p_1)} - c_{\beta, 1}^{(p_2)}\theta_{\beta, 1}^{m_2-m_1}} \theta_{\alpha, 1}^{n_1} \theta_{\beta, 1}^{-m_1}.
    \]
If $\Gamma_{3} =1$, then $\dfrac{c_{\alpha, 1}^{(j_1)} - c_{\alpha, 1}^{(j_2)}}{c_{\beta, 1}^{(p_1)} - c_{\beta, 1}^{(p_2)}\theta_{\beta, 1}^{m_2-m_1}} =\dfrac{\theta_{\beta, 1}^{m_1}}{\theta_{\alpha, 1}^{n_1}}$. Applying Lemma \ref{CPZ}, there exists an effectively computable constant $D_{0}$, such that
\begin{align*}
 D_{0} \max(n_{1}, m_{1}) &\leq h \left( \dfrac{c_{\alpha, 1}^{(j_1)} - c_{\alpha, 1}^{(j_2)}}{c_{\beta, 1}^{(p_1)} - c_{\beta, 1}^{(p_2)}\theta_{\beta, 1}^{m_2-m_1}}\right).
\end{align*}
\noindent Applying \eqref{height_of_sum}, \eqref{height_of_product}, and \eqref{height_of_exponents}, we get 

\begin{align*}
 h \left( \dfrac{c_{\alpha, 1}^{(j_1)} - c_{\alpha, 1}^{(j_2)}}{c_{\beta, 1}^{(p_1)} - c_{\beta, 1}^{(p_2)}\theta_{\beta, 1}^{m_2-m_1}}\right) &\leq h\left(c_{\alpha, 1}^{(j_1)} - c_{\alpha, 1}^{(j_2)}\right) + h\left(c_{\beta, 1}^{(p_1)} - c_{\beta, 1}^{(p_2)}\theta_{\beta, 1}^{m_2-m_1}\right)\\
 &\leq h(c_{\alpha, 1}^{(j_1)}) + h(c_{\alpha, 1}^{(j_2)}) 
+ h(c_{\beta, 1}^{(p_1)})+ h(c_{\beta, 1}^{(p_2)})\\ & \ \ \ \ + (m_1 -m_2) \frac{\log \theta_{\beta, 1}}{2} + 2\log 2.
\end{align*}

\noindent As noted above, 
\begin{equation*}
    m_{1}-m_{2} \leq  C_{10}\log n_{1}. 
\end{equation*}

\noindent Hence
\begin{equation}\label{case(ii)_bounded m_1 2}
    \max(n_1, m_1) \ll \log n_{1}.
\end{equation}

\noindent If $\Gamma_{3}\neq 1$, then we apply Lemma \ref{lfl} with $T=3$, $D=4$,
    \begin{align*}
        \delta_1&=\theta_{\alpha, 1},\ \delta_2=\theta_{\beta, 1},\ \delta_3=\frac{c_{\alpha, 1}^{(j_1)} - c_{\alpha, 1}^{(j_2)}}{c_{\beta, 1}^{(p_1)} - c_{\beta, 1}^{(p_2)}\theta_{\beta, 1}^{m_2-m_1}},\\
        k_1&=n_1,\ k_2=-m_1,\ k_3=1,\\
        A_1'&\ll 1,\ A_2'\ll 1,\ A_3' \ll \log n_1.
    \end{align*}
Using \eqref{LFL1_B}, we can take $B=\max(C_{1},1)n_1$. Then 
    \[
  \log|\Gamma_{3}-1|>-C_{12}\log(e \max(C_{1},1)n_1) \log n_1.  
    \]
Comparing with \eqref{lfl3_rhs}, we obtain
\[ m_{1}< \frac{C_{12}(1+\log \max(C_1, 1)+ \log n_{1}) \log n_1 + \log C_{11}}{\log \theta_{\beta,1}}.\] 
Therefore,
    \begin{equation}\label{case(ii)_bounded m_1 3}
    m_{1} < C_{13} (\log n_1)^{2},
    \end{equation}
where $C_{13}=\dfrac{2 C_{12}\left(2+\log^+ C_1\right)+3\log^{+} C_{11}}{\log \theta_{\beta,1}}$ since $n_1 \geq 2$.

\noindent From inequalities \eqref{case(ii)_bounded m_1 1}, \eqref{case(ii)_bounded m_1 2}, and  \eqref{case(ii)_bounded m_1 3}, we obtain
$$m_1 \ll (\log n_1)^{2} \ll (\log n_1)^{3}.$$

\noindent \textbf{Case (iii):} It is similar to case (ii). Following a similar procedure, we get 
\[ n_1 \ll (\log n_1)^{2}.\]
\noindent Hence, $$m_1 \ll (\log n_1)^{2} \ll (\log n_1)^{3}.$$

\noindent \textbf{Case (iv):} $n_{1} > n_{2}$ and $m_{1} > m_{2}$.\\
\noindent Using equation \eqref{main_eqn_1}, we get
\begin{equation*}
c_{\alpha, 1}^{(j_1)}\theta_{\alpha, 1}^{n_1} - c_{\beta, 1}^{(p_1)} \theta_{\beta, 1}^{m_1} = 
c_{\alpha, 2}^{(j_1)} \theta_{\alpha, 2}^{n_1} + c_{\alpha, 1}^{(j_2)} \theta_{\alpha, 1}^{n_2}- c_{\alpha, 2}^{(j_2)} \theta_{\alpha, 2}^{n_2} -  c_{\beta, 2}^{(p_1)} \theta_{\beta, 2}^{m_1} - c_{\beta, 1}^{(p_2)} \theta_{\beta, 1}^{m_2}+ c_{\beta, 2}^{(p_2)} \theta_{\beta, 2}^{m_2}.    
\end{equation*}
Since $|\theta_{\alpha, 2}|<1$, $|\theta_{\beta, 2}|<1$, we get
\[ \left|c_{\alpha, 1}^{(j_1)}\theta_{\alpha, 1}^{n_1} - c_{\beta, 1}^{(p_1)} \theta_{\beta, 1}^{m_1} \right| \leq  C_{14} +  c_{\alpha, 1}^{(j_2)}  \theta_{\alpha, 1}^{n_2} + c_{\beta, 1}^{(p_2)} \theta_{\beta, 1}^{m_2}.\]
Dividing by $c_{\beta, 1}^{(p_1)}  \theta_{\beta, 1}^{m_1}$ and using \eqref{domterm_samesize}, we get
\begin{align}\nonumber
   \left| \frac{c_{\alpha, 1}^{(j_1)}}{c_{\beta, 1}^{(p_1)}} \theta_{\alpha, 1}^{n_1} \theta_{\beta, 1}^{-m_1}- 1  \right| &\leq 
   \frac{C_{14}}{c_{\beta,1}^{(p_1)}\theta_{\beta, 1}^{m_1}}  + \frac{c_{\alpha, 1}^{(j_2)}}{c_{\beta,1}^{(p_1)}} \frac{\theta_{\alpha, 1}^{n_2}}{\theta_{\beta, 1}^{m_1}} + \frac{c_{\beta, 1}^{(p_2)}}{c_{\beta,1}^{(p_1)}} \frac{\theta_{\beta, 1}^{m_2}}{\theta_{\beta, 1}^{m_1}}\\
  \nonumber 
&\leq \frac{C_{14}}{c_{\beta,1}^{(p_1)}\theta_{\beta, 1}^{m_1}}  + \frac{c_{\alpha, 1}^{(j_2)}}{c_{\beta,1}^{(p_1)}} \frac{\theta_{\alpha, 1}^{n_1}}{\theta_{\beta, 1}^{m_1}} \frac{\theta_{\alpha, 1}^{n_2}}{\theta_{\alpha, 1}^{n_1}}+ \frac{c_{\beta, 1}^{(p_2)}}{c_{\beta,1}^{(p_1)}} \frac{1}{\theta_{\beta, 1}^{m_1-m_2}}\\
\nonumber
&\leq \frac{C_{14}}{c_{\beta,1}^{(p_1)}\theta_{\beta, 1}^{m_1}} + \frac{c_{\alpha, 1}^{(j_2)}}{c_{\beta,1}^{(p_1)}} \frac{C_{\beta, 4}}{C_{\alpha, 3}} \frac{1}{\theta_{\alpha, 1}^{n_{1}-n_{2}}} + \frac{c_{\beta, 1}^{(p_2)}}{c_{\beta,1}^{(p_1)}} \frac{1}{\theta_{\beta, 1}^{m_1-m_2}}\\
\nonumber
&\leq C_{15} \max \left\{\frac{1}{\theta^{n_{1}-n_{2}}_{\alpha,1}}, \frac{1}{\theta^{m_{1}-m_{2}}_{\beta,1 }}\right\} \\
 \label{lfl4_rhs} &\leq \frac{C_{15}}{\min(\theta_{\alpha,1},\theta_{\beta,1})^{\min(n_1-n_2,m_1-m_2)}}.
\end{align}
Let 
    \[
\Gamma_{4}=\frac{c_{\alpha,1}^{(j_1)}\theta_{\alpha,1}^{n_1}}{c_{\beta,1}^{(p_1)}\theta_{\beta,1}^{m_1}}.
    \]
If $\Gamma_{4}=1$, then $\dfrac{c_{\alpha,1}^{(j_{1})}}{ c_{\beta, 1}^{(p_{1})}}=\dfrac{\theta_{\beta, 1}^{m_{1}}}{\theta_{\alpha, 1}^{n_{1}}}$. Taking heights and applying Lemma \ref{CPZ}, we get

\begin{equation}\label{case(iv)_bounded m_1 1}
  \max (n_{1}, m_{1}) \ll 1.  
\end{equation}

\noindent
If $\Gamma_{4}\neq 1$, then we apply Lemma \ref{lfl} with $T=3$, $D=4$,
    \begin{align*}
        \delta_1&=\theta_{\alpha, 1},\ \delta_2=\theta_{\beta, 1},\ \delta_3=\dfrac{c_{\alpha,1}^{(j_{1})}}{ c_{\beta, 1}^{(p_{1})}},\\
        k_1&=n_1,\ k_2=-m_1,\ k_3=1,\\
        A_1'&\ll 1,\ A_2'\ll 1,\ A_3' \ll 1.
    \end{align*}
Using \eqref{LFL1_B}, we can take $B=\max(C_{1},1)n_1$. Then 
    \[
\log|\Gamma_{4}-1|>-C_{16}\log(e \max(C_{1},1)n_1).  
    \]
Comparing with \eqref{lfl4_rhs}, we obtain
\[ \min(n_1-n_2,m_1-m_2)< \frac{C_{16}(1+\log \max(C_1, 1)+ \log n_{1})+\log C_{15}}{\log \min(\theta_{\alpha,1}, \theta_{\beta,1})}.\] 
Therefore,
    \begin{equation*}
    \min(n_1-n_2,m_1-m_2)<C_{17}\log n_1,
    \end{equation*}
where $C_{17}=\dfrac{2 C_{16}\left(2+\log^+ C_1\right)+ 2\log^{+} C_{15}}{\log \min(\theta_{\alpha,1}, \theta_{\beta,1})}$ since $n_1 \geq 2$.\\

\noindent If $\min (n_1-n_2,m_1-m_2) = n_1-n_2,$ then, using equation \eqref{main_eqn_1}, we have
\begin{align}
\nonumber
 c_{\alpha, 1}^{(j_1)}\theta_{\alpha, 1}^{n_1} - c_{\alpha, 1}^{(j_2)}\theta_{\alpha, 1}^{n_2} - c_{\beta, 1}^{(p_1)} \theta_{\beta, 1}^{m_1} 
&=  c_{\alpha, 2}^{(j_1)}\theta_{\alpha, 2}^{n_1} - c_{\alpha, 2}^{(j_2)}\theta_{\alpha, 2}^{n_2}
 - c_{\beta, 2}^{(p_1)} \theta_{\beta, 2}^{m_1} - c_{\beta, 1}^{(p_2)} \theta_{\beta, 1}^{m_2} \\ \nonumber &+ c_{\beta, 2}^{(p_2)} \theta_{\beta, 2}^{m_2}.  
\end{align}
Since $|\theta_{\alpha, 2}|<1$, $|\theta_{\beta, 2}|<1$, we get

\[ \left| c_{\alpha, 1}^{(j_1)}\theta_{\alpha, 1}^{n_1} - c_{\alpha, 1}^{(j_2)}\theta_{\alpha, 1}^{n_2} - c_{\beta, 1}^{(p_1)} \theta_{\beta, 1}^{m_1}  \right| \leq  C_{18} + c_{\beta, 1}^{(p_2)} \theta_{\beta, 1}^{m_2}.\]
Dividing by $c_{\beta, 1}^{(p_1)}  \theta_{\beta, 1}^{m_1}$, we get
\begin{equation}\label{lfl5_rhs}
    \left| \frac{c_{\alpha, 1}^{(j_1)}- c_{\alpha, 1}^{(j_2)}\theta_{\alpha, 1}^{n_2 - n_1}}{c_{\beta, 1}^{(p_1)}} \theta_{\alpha, 1}^{n_1} \theta_{\beta, 1}^{-m_1} -1 \right| \leq  \frac{C_{19}}{\theta_{\beta, 1}^{m_1-m_2}}.
\end{equation}

\noindent Let 
    \[
    \Gamma_{5}=\frac{c_{\alpha, 1}^{(j_1)}- c_{\alpha, 1}^{(j_2)}\theta_{\alpha, 1}^{n_2 - n_1}}{c_{\beta, 1}^{(p_1)}} \theta_{\alpha, 1}^{n_1} \theta_{\beta, 1}^{-m_1}.
    \]
If $\Gamma_{5} =1$, then $\dfrac{c_{\alpha, 1}^{(j_1)}- c_{\alpha, 1}^{(j_2)}\theta_{\alpha, 1}^{n_2 - n_1}}{c_{\beta, 1}^{(p_1)}} =\dfrac{\theta_{\beta, 1}^{m_1}}{\theta_{\alpha, 1}^{n_1}}$. Applying Lemma \ref{CPZ}, there exists an effectively computable constant $D_{0}$, such that 
\begin{align*}
 D_{0} \max(n_{1}, m_{1}) &\leq  h \left(\dfrac{c_{\alpha, 1}^{(j_1)}- c_{\alpha, 1}^{(j_2)}\theta_{\alpha, 1}^{n_2 - n_1}}{c_{\beta, 1}^{(p_1)}}  \right).
\end{align*}

\noindent Applying \eqref{height_of_sum}, \eqref{height_of_product}, and \eqref{height_of_exponents}, we get 
\begin{align*}
 h \left(\dfrac{c_{\alpha, 1}^{(j_1)}- c_{\alpha, 1}^{(j_2)}\theta_{\alpha, 1}^{n_2 - n_1}}{c_{\beta, 1}^{(p_1)}}  \right) &\leq h\left(c_{\alpha, 1}^{(j_1)}- c_{\alpha, 1}^{(j_2)}\theta_{\alpha, 1}^{n_2 - n_1} \right) + h ( c_{\beta, 1}^{(p_1)})\\
 &\leq  h(c_{\alpha, 1}^{(j_1)}) + h(c_{\alpha, 1}^{(j_2)}) + (n_1 -n_2) \frac{\log \theta_{\alpha, 1}}{2} +  h(c_{\beta, 1}^{(p_1)}) \\\ & \ \  \ + \log 2.
\end{align*}


\noindent As noted above, 
\begin{equation*}
    n_{1}-n_{2} \leq C_{17}\log n_{1}. 
\end{equation*}

\noindent Hence
\begin{equation}\label{case(iv)_bounded m_1 2}
    \max(n_1, m_1) \ll \log n_{1}.
\end{equation}

\noindent
If $\Gamma_{5}\neq 1$, then we apply Lemma \ref{lfl} with $T=3$, $D=4$,
    \begin{align*}
        \delta_1&=\theta_{\alpha, 1},\ \delta_2=\theta_{\beta, 1},\ \delta_3=\dfrac{c_{\alpha, 1}^{(j_1)}- c_{\alpha, 1}^{(j_2)}\theta_{\alpha, 1}^{n_2 - n_1}}{c_{\beta, 1}^{(p_1)}},\\
        k_1&=n_1,\ k_2=-m_1,\ k_3=1,\\
        A_1'&\ll 1,\ A_2' \ll 1,\ A_3'\ll \log n_1.
    \end{align*}
Using \eqref{LFL1_B}, we can take $B=\max(C_{1},1)n_1$. Then 
    \[
  \log|\Gamma_{5}-1|>-C_{20}\log(e \max(C_{1},1)n_1) \log n_1.  
    \]
Comparing with \eqref{lfl5_rhs}, we obtain
\[ m_1 - m_2< \frac{C_{20}(1+\log \max(C_1, 1)+ \log n_{1}) \log n_1 +\log C_{19}}{\log  \theta_{\beta,1}}.\] 
Therefore,
    \begin{equation*}
    m_1-m_2 < C_{21} (\log n_1)^{2},
    \end{equation*}
where $C_{21}=\dfrac{2 C_{20}\left(2+\log^+ C_1\right)+ 3\log^{+} C_{19}}{\log \theta_{\beta,1}}$ since $n_1 \geq 2$.\\

\noindent If $\min(n_1 - n_2, m_1 -m_2) = m_1 - m_2,$ then, we get 
    \begin{equation*}
    n_1-n_2 < C_{22} (\log n_1)^{2}.
    \end{equation*}
    
\noindent Using both cases,
\[\max(n_1 - n_2, m_1 -m_2) \ll (\log n_1)^{2}.\]

\noindent By again considering equation \eqref{main_eqn_1}, we get
\begin{align}
\nonumber
 c_{\alpha, 1}^{(j_1)}\theta_{\alpha, 1}^{n_1} - c_{\alpha, 1}^{(j_2)}\theta_{\alpha, 1}^{n_2} - c_{\beta, 1}^{(p_1)} \theta_{\beta, 1}^{m_1} + c_{\beta, 1}^{(p_2)} \theta_{\beta, 1}^{m_2}
&=  c_{\alpha, 2}^{(j_1)}\theta_{\alpha, 2}^{n_1} - c_{\alpha, 2}^{(j_2)}\theta_{\alpha, 2}^{n_2}
 - c_{\beta, 2}^{(p_1)} \theta_{\beta, 2}^{m_1}\\ \nonumber &+ c_{\beta, 2}^{(p_2)} \theta_{\beta, 2}^{m_2}.  
\end{align}
Since $|\theta_{\alpha, 2}|<1$, $|\theta_{\beta, 2}|<1$, we get

\[ \left| c_{\alpha, 1}^{(j_1)}\theta_{\alpha, 1}^{n_1} - c_{\alpha, 1}^{(j_2)}\theta_{\alpha, 1}^{n_2} - c_{\beta, 1}^{(p_1)} \theta_{\beta, 1}^{m_1} + c_{\beta, 1}^{(p_2)} \theta_{\beta, 1}^{m_2} \right| \leq  C_{23}. \]
Dividing by $\left|c_{\beta, 1}^{(p_1)} \theta_{\beta, 1}^{m_1} - c_{\beta, 1}^{(p_2)} \theta_{\beta, 1}^{m_2}\right|
$ and recalling Case (ii), we get
\begin{align}
\nonumber
    \left| \frac{c_{\alpha, 1}^{(j_1)}- c_{\alpha, 1}^{(j_2)}\theta_{\alpha, 1}^{n_2 - n_1}}{c_{\beta, 1}^{(p_1)}  - c_{\beta, 1}^{(p_2)} \theta_{\beta, 1}^{m_2-m_1}} \theta_{\alpha, 1}^{n_1} \theta_{\beta, 1}^{-m_1} -1 \right| &\leq  \frac{C_{23}}{\left| c_{\beta, 1}^{(p_1)}  - c_{\beta, 1}^{(p_2)} \theta_{\beta, 1}^{m_2-m_1}\right|} \frac{1}{\theta_{\beta, 1}^{m_1}}\\
\label{lfl6_rhs}
    &\leq  \frac{C_{24}}{\theta_{\beta, 1}^{m_1}}.
\end{align}

\noindent Let 
    \[
    \Gamma_{6}=\dfrac{c_{\alpha, 1}^{(j_1)}- c_{\alpha, 1}^{(j_2)}\theta_{\alpha, 1}^{n_2 - n_1}}{c_{\beta, 1}^{(p_1)}  - c_{\beta, 1}^{(p_2)} \theta_{\beta, 1}^{m_2-m_1}} \theta_{\alpha, 1}^{n_1} \theta_{\beta, 1}^{-m_1}.
    \]
If $\Gamma_{6} =1$, then $\dfrac{c_{\alpha, 1}^{(j_1)}- c_{\alpha, 1}^{(j_2)}\theta_{\alpha, 1}^{n_2 - n_1}}{c_{\beta, 1}^{(p_1)}  - c_{\beta, 1}^{(p_2)} \theta_{\beta, 1}^{m_2-m_1}} =\dfrac{\theta_{\beta, 1}^{m_1}}{\theta_{\alpha, 1}^{n_1}}$. Applying Lemma \ref{CPZ}, there exists an effectively computable constant $D_{0}$, such that 
\begin{align*}
 D_{0} \max(n_{1}, m_{1}) &\leq  h \left(\frac{c_{\alpha, 1}^{(j_1)}- c_{\alpha, 1}^{(j_2)}\theta_{\alpha, 1}^{n_2 - n_1}}{c_{\beta, 1}^{(p_1)}  - c_{\beta, 1}^{(p_2)} \theta_{\beta, 1}^{m_2-m_1}}  \right).
\end{align*}
\noindent Applying \eqref{height_of_sum}, \eqref{height_of_product}, and \eqref{height_of_exponents}, we get 

\begin{align*}
 h \left(\frac{c_{\alpha, 1}^{(j_1)}- c_{\alpha, 1}^{(j_2)}\theta_{\alpha, 1}^{n_2 - n_1}}{c_{\beta, 1}^{(p_1)}  - c_{\beta, 1}^{(p_2)} \theta_{\beta, 1}^{m_2-m_1}}  \right) &\leq h\left(c_{\alpha, 1}^{(j_1)}- c_{\alpha, 1}^{(j_2)}\theta_{\alpha, 1}^{n_2 - n_1} \right) + h\left(c_{\beta, 1}^{(p_1)}  - c_{\beta, 1}^{(p_2)} \theta_{\beta, 1}^{m_2-m_1}\right)\\
 &\leq h(c_{\alpha, 1}^{(j_1)}) + h(c_{\alpha, 1}^{(j_2)}) 
+ h(c_{\beta, 1}^{(p_1)})+ h(c_{\beta, 1}^{(p_2)}) + 2\log 2  \\ & \ \ \ + (m_1 -m_2) \frac{\log \theta_{\beta, 1}}{2} + (n_1 -n_2) \frac{\log \theta_{\alpha, 1}}{2}.
\end{align*}



\noindent As noted above, 
\begin{equation*}
    \max (n_{1}-n_{2}, m_{1}-m_{2}) \ll (\log n_{1})^{2}. 
\end{equation*}

\noindent Hence
\begin{equation}\label{case(iv)_bounded m_1 3}
    \max(n_1, m_1) \ll (\log n_{1})^{2}.
\end{equation}

\noindent
If $\Gamma_{6}\neq 1$, then we apply Lemma \ref{lfl} with $T=3$, $D=4$,
    \begin{align*}
        \delta_1&=\theta_{\alpha, 1},\ \delta_2=\theta_{\beta, 1},\ \delta_3=\frac{c_{\alpha, 1}^{(j_1)}- c_{\alpha, 1}^{(j_2)}\theta_{\alpha, 1}^{n_2 - n_1}}{c_{\beta, 1}^{(p_1)}  - c_{\beta, 1}^{(p_2)} \theta_{\beta, 1}^{m_2-m_1}},\\
        k_1&=n_1,\ k_2=-m_1,\ k_3=1,\\
        A_1'& \ll 1,\ A_2' \ll 1,\ A_3' \ll (\log n_1)^{2}.
    \end{align*}
Using \eqref{LFL1_B}, we can take $B=\max(C_{1},1)n_1$. Then 
    \[
  \log|\Gamma_{6}-1|>-C_{26}\log(e \max(C_{1},1)n_1) (\log n_1)^{2}.  
    \]
Comparing with \eqref{lfl6_rhs}, we obtain
\[ m_1 < \frac{C_{26}(1+\log \max(C_1, 1)+ \log n_{1}) (\log n_1)^{2}+\log C_{24}}{\log \theta_{\beta,1}}.\] 
Therefore,
    \begin{equation}\label{case(iv)_bounded m_1 4}
    m_1 < C_{27} (\log n_1)^{3},
    \end{equation}
where $C_{27}=\dfrac{2 C_{26}\left(2+\log^+ C_1 \right)+ 4\log^{+} C_{24}}{\log \theta_{\beta,1}}$ since $n_1 \geq 2$.\\ 

\noindent From inequalities \eqref{case(iv)_bounded m_1 1}, \eqref{case(iv)_bounded m_1 2}, \eqref{case(iv)_bounded m_1 3}, and \eqref{case(iv)_bounded m_1 4}, we obtain 
\[m_1 \ll (\log n_1)^{3}.\]

\noindent Now, using equation \eqref{domterm_samesize}, we get
\begin{equation*}
    n_1 < C_{28} (\log n_1)^{3}.
\end{equation*}

\noindent Therefore,  $n_1 <  (C_{29} \log n_1)^{3}.$
    We now use Lemma \ref{PW} with $a=0, c=3, g=C_{29}^{3}= C_{28}$. (If necessary, we can replace $C_{29}$ by $e^2$.) We thus obtain
\begin{align*}
   n_{1} &\leq  216 C_{29}^{3} (\log( 3 C_{29}))^{3}.
\end{align*}

   
\noindent Since $n_1$ is bounded, $m_1$ will also be bounded. Hence, $N_1$ and $M_1$ are bounded. Therefore, the cardinality of the set $\mathcal{C}$ is finite. Hence, Theorem \ref{main thm} is proved.

\section{Proof of Theorem \ref{thm on bound}}\label{Sec_thm on bound} 
Let $c$ be an integer which has two distinct representations as shown in equation \eqref{dist_rep_c}. We note that for our choices of $\alpha$ and $\beta$, we have $s_0 = s_1 =2$ and $r_0 = r_1 = 0$, thus, $N_{i} = 2n_{i} + j_i,$ $ 0 \leq j_i \leq 1$ and $M_i = 2 m_i + p_i, \ 0 \leq p_i \leq 1, \  i=1, 2.$ For technical reasons, we assume that $n_{1} \geq 6$. 
\noindent Since $\alpha = [0;\overline{1, a}]$, we have $q_{\alpha, 0} =1=q_{\alpha, 1}$ and 
\[q_{\alpha, N}= 
\begin{cases*}
a q_{\alpha, N-1}+ q_{\alpha, N-2} & \text{if}\; $N$ is even,\\
q_{\alpha, N-1}+ q_{\alpha, N-2} & \text{if}\; $N$ is odd. 
\end{cases*}\]
Let $d_{a} = a^2 + 4a$ and $\theta_{\alpha, 1} = \dfrac{a+2+\sqrt{d_{a}}}{2}$. Then
\begin{align*}
    q_{\alpha, 2\ell} &= \frac{a+\sqrt{d_{a}}}{2 \sqrt{d_{a}}} \theta_{\alpha,1}^{\ell} - \frac{a-\sqrt{d_{a}}}{2 \sqrt{d_{a}}} \theta_{\alpha,1}^{-\ell},\\
    q_{\alpha, 2\ell+1} &= \frac{\theta_{\alpha, 1}}{\sqrt{d_{a}}} \theta_{\alpha,1}^{\ell} - \frac{\theta_{\alpha, 1}^{-1}}{\sqrt{d_{a}}} \theta_{\alpha,1}^{-\ell}.
\end{align*}
We have the following inequalities:
\begin{align}
\label{q_alpha_eqn1}
    \frac{1}{2} \theta_{\alpha, 1}^{\ell} &< q_{\alpha, 2\ell} < \frac{11}{10} \frac{a+\sqrt{d_{a}}}{2\sqrt{d_{a}}} \theta_{\alpha, 1}^{\ell},\\
\label{q_alpha_eqn2}
    \frac{7}{8} \frac{\theta_{\alpha,1}}{\sqrt{d_{a}}} \theta_{\alpha, 1}^{\ell} &< q_{\alpha, 2\ell +1} < \frac{9}{8} \frac{\theta_{\alpha,1}}{\sqrt{d_{a}}} \theta_{\alpha, 1}^{\ell}.
\end{align}

\noindent From inequalities \eqref{q_alpha_eqn1} and \ref{q_alpha_eqn2}, we have
\begin{equation*}
    \frac{1}{2} \theta_{\alpha, 1}^{\lfloor N/2 \rfloor} < q_{\alpha, N} < \frac{9}{8} \frac{\theta_{\alpha,1}}{\sqrt{d_{a}}} \theta_{\alpha, 1}^{\lfloor N/2 \rfloor}, \ N \geq 0.
\end{equation*}
We note that $\dfrac{\frac{9 \theta_{\alpha,1}}{8 \sqrt{d_{a}}} }{\nicefrac{1}{2}} < \theta_{\alpha,1}$.
Referring to equation\eqref{lb_3}, 
\begin{align*}
    \frac{1}{2 \sqrt{d_{a}}} \theta_{\alpha, 1}^{\lfloor N_{1}/2 \rfloor} \leq |q_{\alpha, N_1} - q_{\alpha, N_2}| = |q_{\beta, M_1} - q_{\beta, M_2}| < q_{\beta, M_1} 
    < \frac{9}{8} \frac{\theta_{\beta,1}}{\sqrt{d_{b}}} \theta_{\beta, 1}^{\lfloor M_{1}/2 \rfloor}.
\end{align*}     
As $\lfloor \frac{N_1}{2} \rfloor = n_1$ and $\lfloor \frac{M_1}{2} \rfloor = m_1$, we have

\begin{equation}\label{relation_n1m1}
    \frac{1}{2 \sqrt{d_{a}}} \theta_{\alpha, 1}^{n_1} < \frac{9}{8} \frac{\theta_{\beta,1}}{\sqrt{d_{b}}} \theta_{\beta, 1}^{m_1},
\end{equation} 

\noindent therefore
\[ n_1 < \frac{\log(\frac{9}{8} \frac{\theta_{\beta,1}}{\sqrt{d_{b}}}) + \log(2\sqrt{d_{a}})}{\log \theta_{\alpha,1}} + m_1 \frac{\log \theta_{\beta,1}}{\log \theta_{\alpha,1}}.\]
As $y = \log\left(\dfrac{9}{8}\dfrac{x+2+\sqrt{x^2+4x}}{2\sqrt{x^2+4x}}\right)$ and $y = \frac{0.17 + \log(2\sqrt{x^2 + 4x})}{\log (\frac{x+2+\sqrt{x^2+4x}}{2})}$ are decreasing functions of $x$ for $x \geq 2$, we get $\log\left(\dfrac{9}{8} \frac{\theta_{\beta,1}}{\sqrt{d_{b}}}\right) < 0.17$ for $b \geq 3$, and $\frac{0.17 + \log(2\sqrt{d_{a}})}{\log \theta_{\alpha,1}} < 1.60$ for $a \geq 2$.
Therefore,
\begin{equation}\label{n1Rm1_eqn1}
    n_1 < 1.60 + m_1 \frac{\log \theta_{\beta,1}}{\log \theta_{\alpha,1}}.
\end{equation} 
Similarly,
\[ m_1 < \frac{\log(\frac{9}{8} \frac{\theta_{\alpha,1}}{\sqrt{d_{a}}}) + 
\log(2\sqrt{d_{b}})}{\log \theta_{\beta,1}} + n_1 \frac{\log \theta_{\alpha,1}}{\log \theta_{\beta,1}}.\]
Since $\frac{\log(\frac{9}{8} \frac{\theta_{\alpha,1}}{\sqrt{d_{a}}}) + 
\log(2\sqrt{d_{b}})}{\log \theta_{\beta,1}} < 1.55$ for $2 \leq a <b$, we have
\begin{equation}\label{m1Rn1_eqn2}
   m_1 < 1.55 + n_1 \frac{\log \theta_{\alpha,1}}{\log \theta_{\beta,1}}. 
\end{equation}

As in the proof of Theorem \ref{main thm}, the bound obtained in Case (iv) is the largest among all cases; therefore, we discuss only Case (iv).\\
\textbf{Case (iv):} $n_1 > n_2$ and $m_1 > m_2$.\\
\noindent Using equation \eqref{main_eqn_1}, we get
\begin{align*}
    \left|c_{\alpha, 1}^{(j_1)}\theta_{\alpha, 1}^{n_1} - c_{\beta, 1}^{(p_1)} \theta_{\beta, 1}^{m_1} \right| \leq \left| c_{\alpha, 2}^{(j_1)}\right| + \left|c_{\alpha, 2}^{(j_2)}\right|+ \left|c_{\beta, 2}^{(p_1)}\right|+ \left|c_{\beta, 2}^{(p_2)}\right|+ c_{\alpha, 1}^{(j_2)}  \theta_{\alpha, 1}^{n_2} + c_{\beta, 1}^{(p_2)} \theta_{\beta, 1}^{m_2}.
\end{align*}
Since $\dfrac{\sqrt{d_{a}}-a}{2\sqrt{d_{a}}} > \dfrac{\theta_{\alpha,1}^{-1}}{\sqrt{d_{a}}}$, $y = \dfrac{\sqrt{x^2 + 4x} - x}{2 \sqrt{x^2 +4x}}$ is a decreasing function of $x$ for $x \geq 2$, and $2 \leq a < b$, we get \begin{align}
\label{ubof4coeffs} 
\left| c_{\alpha, 2}^{(j_1)}\right| + \left|c_{\alpha, 2}^{(j_2)}\right|+ \left|c_{\beta, 2}^{(p_1)}\right|+ \left|c_{\beta, 2}^{(p_2)}\right| \leq 2  \left(\frac{\sqrt{d_{a}}-a}{2\sqrt{d_{a}}} +   \frac{\sqrt{d_{b}}-b}{2\sqrt{d_{b}}}\right) \leq 0.77.
\end{align}
Therefore,
\begin{align*}
    \left|c_{\alpha, 1}^{(j_1)}\theta_{\alpha, 1}^{n_1} - c_{\beta, 1}^{(p_1)} \theta_{\beta, 1}^{m_1} \right| \leq 0.77+ c_{\alpha, 1}^{(j_2)}  \theta_{\alpha, 1}^{n_2} + c_{\beta, 1}^{(p_2)} \theta_{\beta, 1}^{m_2}.
\end{align*}
Dividing by $c_{\beta, 1}^{(p_1)} \theta_{\beta, 1}^{m_1}$, we get
\begin{align*}
    \left|\frac{c_{\alpha, 1}^{(j_1)}\theta_{\alpha, 1}^{n_1}}{c_{\beta, 1}^{(p_1)} \theta_{\beta, 1}^{m_1}} - 1 \right| \leq \frac{0.77}{c_{\beta, 1}^{(p_1)} \theta_{\beta, 1}^{m_1}} + \frac{c_{\alpha, 1}^{(j_2)}  \theta_{\alpha, 1}^{n_2}}{c_{\beta, 1}^{(p_1)} \theta_{\beta, 1}^{m_1}}  + \frac{c_{\beta, 1}^{(p_2)} \theta_{\beta, 1}^{m_2}}{c_{\beta, 1}^{(p_1)} \theta_{\beta, 1}^{m_1}}.
\end{align*}
\noindent First, we compute an upper bound for the term $\dfrac{0.77}{c_{\beta, 1}^{(p_1)} \theta_{\beta, 1}^{m_1}} + \dfrac{c_{\beta, 1}^{(p_2)} \theta_{\beta, 1}^{m_2}}{c_{\beta, 1}^{(p_1)} \theta_{\beta, 1}^{m_1}}$ by using that the functions $y= \dfrac{\sqrt{x^2 + 4x}}{x+\sqrt{x^2 + 4x}}$ and $y= \dfrac{x+2+\sqrt{x^2 + 4x}}{x+\sqrt{x^2 + 4x}}$ are decreasing functions of $x$ for $x \geq 3$. So,
\begin{align}
\nonumber
  \frac{0.77}{c_{\beta, 1}^{(p_1)} \theta_{\beta, 1}^{m_1}} + \frac{c_{\beta, 1}^{(p_2)} \theta_{\beta, 1}^{m_2}}{c_{\beta, 1}^{(p_1)} \theta_{\beta, 1}^{m_1}} &\leq  \frac{0.77}{\frac{b+\sqrt{d_{b}}}{2\sqrt{d_{b}}}} \frac{1}{\theta_{\beta, 1}^{m_{1}}}+ \frac{\frac{\theta_{\beta, 1}}{\sqrt{d_{b}}}}{\frac{b+\sqrt{d_{b}}}{2\sqrt{d_{b}}}} \frac{1}{\theta_{\beta, 1}^{m_{1}-m_{2}}}\\
  \nonumber
  &=\frac{1.54 \sqrt{d_{b}}}{b+\sqrt{d_{b}}} \frac{1}{\theta_{\beta, 1}^{m_{1}}}+ \frac{b+2+\sqrt{d_{b}}}{b+\sqrt{d_{b}}} \frac{1}{\theta_{\beta, 1}^{m_{1}-m_{2}}}\\
  \nonumber
  &\leq \frac{0.94}{\theta_{\beta, 1}^{m_{1}}}+ \frac{1.27}{\theta_{\beta, 1}^{m_{1}-m_{2}}}\\
  \label{ub1stterm}   &\leq \frac{2.21}{\theta_{\beta, 1}^{m_{1}-m_{2}}}  \ \ \ (\text{because} \ \ m_1 > m_2). 
\end{align}


\noindent Next, we compute an upper bound for $\dfrac{c_{\alpha, 1}^{(j_2)}  \theta_{\alpha, 1}^{n_2}}{c_{\beta, 1}^{(p_1)} \theta_{\beta, 1}^{m_1}}$ by using that the function $y=\dfrac{x+2+\sqrt{x^2+4x}}{2(x+\sqrt{x^2+4x})}$ is a decreasing function of $x$ for $x\geq 3$. So,
\begin{align}
\nonumber
    \frac{c_{\alpha, 1}^{(j_2)}  \theta_{\alpha, 1}^{n_2}}{c_{\beta, 1}^{(p_1)} \theta_{\beta, 1}^{m_1}} =  \frac{c_{\alpha, 1}^{(j_2)}} {c_{\beta, 1}^{(p_1)}}  \frac{\theta_{\alpha, 1}^{n_1}}{\theta_{\beta, 1}^{m_1}} \frac{\theta_{\alpha, 1}^{n_2}}{\theta_{\alpha, 1}^{n_1}}  &\leq \frac{\theta_{\alpha,1}}{\sqrt{d_{a}}} \frac{2 \sqrt{d_{b}}}{b+\sqrt{d_{b}}} \frac{\theta_{\alpha, 1}^{n_1}}{\theta_{\beta, 1}^{m_1}} \frac{1}{\theta_{\alpha, 1}^{n_{1} - n_{2}}}\\
\nonumber
    &\leq \frac{\theta_{\alpha,1}}{\sqrt{d_{a}}} \frac{2 \sqrt{d_{b}}}{b+\sqrt{d_{b}}} \frac{9 \theta_{\beta,1} \sqrt{d_{a}}}{4 \sqrt{d_{b}} } \frac{1}{\theta_{\alpha, 1}^{n_{1} - n_{2}}} \ \ \ (\text{using equation} \ \ \eqref{relation_n1m1}\\
\label{ub2ndTerm}
    &= \frac{9 \theta_{\alpha,1} \theta_{\beta, 1}}{2(b+\sqrt{d_{b}})} \frac{1}{\theta_{\alpha, 1}^{n_{1} - n_{2}}} \leq  \frac{2.85 \  \theta_{\alpha,1}}{\theta_{\alpha,1}^{n_{1}-n_{2}}}.
\end{align}


\noindent Finally, 
\begin{align}
 \nonumber   \left|\frac{c_{\alpha, 1}^{(j_1)}\theta_{\alpha, 1}^{n_1}}{c_{\beta, 1}^{(p_1)} \theta_{\beta, 1}^{m_1}} - 1 \right| &\leq \frac{2.21}{\theta_{\beta, 1}^{m_{1}-m_{2}}}+ \frac{2.85 \theta_{\alpha,1}}{\theta_{\alpha,1}^{n_{1}-n_{2}}}\\
\nonumber &\leq 5.06 \theta_{\alpha,1} \left( \frac{1}{\theta_{\beta, 1}^{m_{1}-m_{2}}}+ \frac{1}{\theta_{\alpha,1}^{n_{1}-n_{2}}} \right)\\
\label{lfl7_rhs}  &\leq 10.12 \theta_{\alpha,1} \max \left( \frac{1}{\theta_{\alpha,1}^{n_{1}-n_{2}}}, \frac{1}{\theta_{\beta, 1}^{m_{1}-m_{2}}} \right).
\end{align}



\noindent
Let $\Gamma_{7} = \dfrac{c_{\alpha, 1}^{(j_1)}\theta_{\alpha, 1}^{n_1}}{c_{\beta, 1}^{(p_1)} \theta_{\beta, 1}^{m_1}}$. If $\Gamma_{7}=1$, then $c_{\alpha, 1}^{(j_1)}\theta_{\alpha, 1}^{n_1}=c_{\beta, 1}^{(p_1)} \theta_{\beta, 1}^{m_1}$ which is not possible because $\mathbb{Q}(\alpha) \neq \mathbb{Q}(\beta)$.
As $\Gamma_{7} \neq 1$, we apply Lemma \ref{lfl} with $T=3$, $D=4$,
    \begin{align*}
        \delta_1&=\theta_{\alpha, 1},\ \delta_2=\theta_{\beta, 1},\ \delta_3=\dfrac{c_{\alpha,1}^{(j_{1})}}{ c_{\beta, 1}^{(p_{1})}},\\
        k_1&=n_1,\ k_2=-m_1,\ k_3=1,\\
        A_1'&=2\log{\theta_{\alpha, 1}},\ A_2'=2\log{\theta_{\beta, 1}}.
    \end{align*}
    Using equation \eqref{height_of_product}, we get
    $h\left(\delta_{3}\right) \leq h(c_{\alpha,1}^{(j_{1})}) + h( c_{\beta, 1}^{(p_{1})})$ .\\
    Now,
    \begin{align*}
     h \left(\frac{a+\sqrt{d_{a}}}{2\sqrt{d_{a}}}\right) &\leq h(a+\sqrt{d_{a}})+h(2\sqrt{d_{a}})
     = \frac{\log 4a}{2} + \frac{\log 4d_{a}}{2},
    \end{align*}
    and
    \begin{align*}
     h \left(\frac{\theta_{\alpha,1}}{\sqrt{d_{a}}}\right) &\leq h(\theta_{\alpha,1})+h(\sqrt{d_{a}}) = \frac{\log \theta_{\alpha,1}}{2} + \frac{\log d_{a}}{2}.
    \end{align*}
    Since $\log 4a + \log 4d_{a} > \log\theta_{\alpha,1} + \log d_{a}$, we have
    $h\left(c_{\alpha, 1}^{(j_1)}\right) \leq \dfrac{\log 4a}{2} + \dfrac{\log 4d_{a}}{2}$. Similarly $h\left(c_{\beta, 1}^{(p_1)}\right) \leq \dfrac{\log 4b}{2} + \dfrac{\log 4d_{b}}{2}$. Since $2 \leq a < b$,  $h\left(\dfrac{c_{\alpha, 1}^{(j_1)}}{c_{\beta, 1}^{(p_1)}}\right) \leq \log 4b  + \log 4d_{b} < 5 \log \theta_{\beta, 1}$. Therefore, we may take $A_3'= 20 \log{\theta_{\beta, 1}}.$ By using equation \eqref{m1Rn1_eqn2}, we can take $B=2n_1$. Then 
    \[
  \log|\Gamma_{7}-1|>- 4.3736 \times 10^{14}(1+\log(2n_{1}))(\log{\theta_{\alpha, 1}}) (\log{\theta_{\beta, 1}})^{2}.  
    \]
Comparing with \eqref{lfl7_rhs}, we obtain
\[\min(n_1-n_2,m_1-m_2) < 8.75 \times 10^{14} (\log{\theta_{\beta, 1}})^{2} \log n_1.\]
\noindent If $\min (n_1-n_2,m_1-m_2) = n_1-n_2,$ then using equation \eqref{main_eqn_1}, we have
\begin{align}
\nonumber
 c_{\alpha, 1}^{(j_1)}\theta_{\alpha, 1}^{n_1} - c_{\alpha, 1}^{(j_2)}\theta_{\alpha, 1}^{n_2} - c_{\beta, 1}^{(p_1)} \theta_{\beta, 1}^{m_1} 
&=  c_{\alpha, 2}^{(j_1)}\theta_{\alpha, 2}^{n_1} - c_{\alpha, 2}^{(j_2)}\theta_{\alpha, 2}^{n_2}
 - c_{\beta, 2}^{(p_1)} \theta_{\beta, 2}^{m_1} - c_{\beta, 1}^{(p_2)} \theta_{\beta, 1}^{m_2} \\ \nonumber &+ c_{\beta, 2}^{(p_2)} \theta_{\beta, 2}^{m_2}. 
\end{align}
Since $|\theta_{\alpha, 2}|<1$, $|\theta_{\beta, 2}|<1$, we get

\begin{align*}
    \left| c_{\alpha, 1}^{(j_1)}\theta_{\alpha, 1}^{n_1} - c_{\alpha, 1}^{(j_2)}\theta_{\alpha, 1}^{n_2} - c_{\beta, 1}^{(p_1)} \theta_{\beta, 1}^{m_1}  \right| &\leq  \left| c_{\alpha, 2}^{(j_1)}\right| + \left|c_{\alpha, 2}^{(j_2)}\right|+ \left|c_{\beta, 2}^{(p_1)}\right|+ \left|c_{\beta, 2}^{(p_2)}\right| + c_{\beta, 1}^{(p_2)} \theta_{\beta, 1}^{m_2}\\
    &\leq 0.77 + c_{\beta, 1}^{(p_2)} \theta_{\beta, 1}^{m_2} \ \ \text{ (using inequality \eqref{ubof4coeffs}).}
\end{align*}
Dividing by $c_{\beta, 1}^{(p_1)} \theta_{\beta, 1}^{m_1}$ and using inequality \eqref{ub1stterm}, we get
\begin{align}
\nonumber
  \left| \frac{c_{\alpha, 1}^{(j_1)}\theta_{\alpha, 1}^{n_1} - c_{\alpha, 1}^{(j_2)}\theta_{\alpha, 1}^{n_2}}{c_{\beta, 1}^{(p_1)} \theta_{\beta, 1}^{m_1}} - 1  \right| &\leq \frac{0.77}{c_{\beta, 1}^{(p_1)} \theta_{\beta, 1}^{m_1}}+ \frac{c_{\beta, 1}^{(p_2)} \theta_{\beta, 1}^{m_2}}{c_{\beta, 1}^{(p_1)} \theta_{\beta, 1}^{m_1}}\\
\label{lfl8_rhs}
  &\leq \frac{2.21}{\theta_{\beta, 1}^{m_1-m_2}}.
\end{align}

\noindent Let 
    \[
    \Gamma_{8}=\frac{c_{\alpha, 1}^{(j_1)}- c_{\alpha, 1}^{(j_2)}\theta_{\alpha, 1}^{n_2 - n_1}}{c_{\beta, 1}^{(p_1)}} \theta_{\alpha, 1}^{n_1} \theta_{\beta, 1}^{-m_1}.
    \]
If $\Gamma_{8} = 1$, then $c_{\alpha, 1}^{(j_1)} \theta_{\alpha, 1}^{n_1} - c_{\alpha, 1}^{(j_2)}\theta_{\alpha, 1}^{n_2} = c_{\beta, 1}^{(p_1)} \theta_{\beta, 1}^{m_1}$ which is not possible because $\mathbb{Q}(\alpha) \neq \mathbb{Q}(\beta).$ As $\Gamma_{8}\neq 1$, we apply Lemma \ref{lfl} with $T=3$, $D=4$,
    \begin{align*}
        \delta_1&=\theta_{\alpha, 1},\ \delta_2=\theta_{\beta, 1},\ \delta_3=\dfrac{c_{\alpha, 1}^{(j_1)}- c_{\alpha, 1}^{(j_2)}\theta_{\alpha, 1}^{n_2 - n_1}}{c_{\beta, 1}^{(p_1)}},\\
        k_1&=n_1,\ k_2=-m_1,\ k_3=1,\\
        A_1'&=2\log{\theta_{\alpha, 1}},\ A_2'=2\log{\theta_{\beta, 1}}.    \end{align*}
    Applying \eqref{height_of_sum}, \eqref{height_of_product}, and \eqref{height_of_exponents}, we get

\begin{align*}
 h \left(\delta_{3}  \right) &\leq h\left(c_{\alpha, 1}^{(j_1)}- c_{\alpha, 1}^{(j_2)}\theta_{\alpha, 1}^{n_2 - n_1} \right) + h\left( c_{\beta, 1}^{(p_1)}\right)\\
 &\leq h(c_{\alpha, 1}^{(j_1)}) + h(c_{\alpha, 1}^{(j_2)}) + h(c_{\beta, 1}^{(p_1)}) + (n_1 - n_2 ) \frac{\log \theta_{\alpha, 1}}{2}+ \log 2. 
\end{align*}


\noindent As
\begin{align*}
  h(c_{\alpha, 1}^{(j_1)}) + h(c_{\alpha, 1}^{(j_2)}) + h(c_{\beta, 1}^{(p_1)}) &\leq \frac{3}{2}(\log 4b + \log 4d_{b}) \leq \frac{15}{2} \log \theta_{\beta, 1} ,\\ \text{ and } 
   n_1 - n_2 &< 8.75 \times 10^{14} (\log{\theta_{\beta, 1}})^{2} \log n_1,
\end{align*}
we get 
\begin{align*}
 h \left(\delta_{3}  \right) & \leq 8.75 \times 10^{14} (\log{\theta_{\beta, 1}})^{2} \log n_1 \frac{\log \theta_{\alpha, 1} }{2}  + \frac{15}{2} \log \theta_{\beta,1} + \log 2 \\
 &< 8.75 \times 10^{14} (\log{\theta_{\beta, 1}})^{3} \log n_1.
\end{align*}
 Therefore, we may take $A_{3}' = 3.5 \times 10^{15} (\log{\theta_{\beta, 1}})^{3} \log n_1.$\\
\noindent Using equation \eqref{m1Rn1_eqn2}, we can take $B=2n_1$. Then 
    \[
  \log|\Gamma_{8}-1|>- 7.3567 \times10^{28}(1+\log 2n_1) \log \theta_{\alpha,1} (\log \theta_{\beta,1})^{4} \log n_1.  
    \]
Comparing with \eqref{lfl8_rhs}, we obtain
\[ m_1 - m_2< 1.54 \times10^{29}  (\log \theta_{\beta,1})^{4} (\log n_1)^2.\] 

\noindent In a similar manner, we can show that if the $\min(n_1 - n_2, m_1 -m_2) = m_1 - m_2$, then $n_1 -n_2 < 1.54 \times10^{29}  (\log \theta_{\beta,1})^{4} (\log n_1)^2.$\\ 
Hence
\[\max(n_1 -n_2, m_1 -m_2) < 1.54 \times10^{29}  (\log \theta_{\beta,1})^{4} (\log n_1)^2.\]

\noindent By again considering equation \eqref{main_eqn_1}, we get
\begin{align}
\nonumber
 c_{\alpha, 1}^{(j_1)}\theta_{\alpha, 1}^{n_1} - c_{\alpha, 1}^{(j_2)}\theta_{\alpha, 1}^{n_2} - c_{\beta, 1}^{(p_1)} \theta_{\beta, 1}^{m_1} + c_{\beta, 1}^{(p_2)} \theta_{\beta, 1}^{m_2}
&=  c_{\alpha, 2}^{(j_1)}\theta_{\alpha, 2}^{n_1} - c_{\alpha, 2}^{(j_2)}\theta_{\alpha, 2}^{n_2}
 - c_{\beta, 2}^{(p_1)} \theta_{\beta, 2}^{m_1}\\ \nonumber &+ c_{\beta, 2}^{(p_2)} \theta_{\beta, 2}^{m_2}.   
\end{align}
Since $|\theta_{\alpha, 2}|<1$, $|\theta_{\beta, 2}|<1$, we get

\begin{align}
 \nonumber   \left| c_{\alpha, 1}^{(j_1)}\theta_{\alpha, 1}^{n_1} - c_{\alpha, 1}^{(j_2)}\theta_{\alpha, 1}^{n_2} - c_{\beta, 1}^{(p_1)} \theta_{\beta, 1}^{m_1} + c_{\beta, 1}^{(p_2)} \theta_{\beta, 1}^{m_2} \right| &\leq  \left| c_{\alpha, 2}^{(j_1)}\right| + \left|c_{\alpha, 2}^{(j_2)}\right|+ \left|c_{\beta, 2}^{(p_1)}\right|+ \left|c_{\beta, 2}^{(p_2)}\right|\\
\label{3rd Ineq}    &\leq 0.77  \ \ \text{ (using inequality \eqref{ubof4coeffs}).}
\end{align}
Dividing by $\left|c_{\beta, 1}^{(p_1)} \theta_{\beta, 1}^{m_1} - c_{\beta, 1}^{(p_2)} \theta_{\beta, 1}^{m_2}\right|
$, we get
\begin{align*}
    \left| \frac{c_{\alpha, 1}^{(j_1)}- c_{\alpha, 1}^{(j_2)}\theta_{\alpha, 1}^{n_2 - n_1}}{c_{\beta, 1}^{(p_1)}  - c_{\beta, 1}^{(p_2)} \theta_{\beta, 1}^{m_2-m_1}} \theta_{\alpha, 1}^{n_1} \theta_{\beta, 1}^{-m_1} -1 \right| \leq  \frac{0.77}{\left| c_{\beta, 1}^{(p_1)}  - c_{\beta, 1}^{(p_2)} \theta_{\beta, 1}^{m_2-m_1}\right|} \frac{1}{\theta_{\beta, 1}^{m_1}}.
\end{align*}
Since $ \left| \frac{1}{\theta_{\beta,1}^{k}}\right| \leq \frac{c_{\beta,1}^{(p_1)}}{2c_{\beta,1}^{(p_2)}}$ holds for all $k \geq 1$ and $y = \frac{x+\sqrt{x^2 + 4x}}{2 \sqrt{x^2+ 4x}}$ is a decreasing function of $x$ for $x \geq 3$, we get 
\begin{align}
\nonumber
    \left| \frac{c_{\alpha, 1}^{(j_1)}- c_{\alpha, 1}^{(j_2)}\theta_{\alpha, 1}^{n_2 - n_1}}{c_{\beta, 1}^{(p_1)}  - c_{\beta, 1}^{(p_2)} \theta_{\beta, 1}^{m_2-m_1}} \theta_{\alpha, 1}^{n_1} \theta_{\beta, 1}^{-m_1} -1 \right| \leq  \frac{0.77 }{\left| \frac{c_{\beta, 1}^{(p_1)} }{2} \right|} \frac{1}{\theta_{\beta, 1}^{m_1}}
    &\leq \frac{1.54}{\frac{b+\sqrt{d_{b}}}{2 \sqrt{d_{b}}}} \frac{1}{\theta_{\beta, 1}^{m_1}}\\
\label{lfl9_rhs} &\leq \frac{1.87}{\theta_{\beta, 1}^{m_1}}.
\end{align}
    
\noindent Let 
    \[
    \Gamma_{9}=\dfrac{c_{\alpha, 1}^{(j_1)}- c_{\alpha, 1}^{(j_2)}\theta_{\alpha, 1}^{n_2 - n_1}}{c_{\beta, 1}^{(p_1)}  - c_{\beta, 1}^{(p_2)} \theta_{\beta, 1}^{m_2-m_1}} \theta_{\alpha, 1}^{n_1} \theta_{\beta, 1}^{-m_1}.
    \]
If $\Gamma_{9} =1$, then $c_{\alpha, 1}^{(j_1)} \theta_{\alpha, 1}^{ n_1} - c_{\alpha, 1}^{(j_2)} \theta_{\alpha, 1}^{n_2} = c_{\beta, 1}^{(p_1)} \theta_{\beta, 1}^{m_1}  - c_{\beta, 1}^{(p_2)} \theta_{\beta, 1}^{m_2}$ which is not possible because $\mathbb{Q}(\alpha) \neq \mathbb{Q}(\beta)$. As $\Gamma_{9}\neq 1$, we apply Lemma \ref{lfl} with $T=3$, $D=4$,
    \begin{align*}
        \delta_1&=\theta_{\alpha, 1},\ \delta_2=\theta_{\beta, 1},\ \delta_3=\frac{c_{\alpha, 1}^{(j_1)}- c_{\alpha, 1}^{(j_2)}\theta_{\alpha, 1}^{n_2 - n_1}}{c_{\beta, 1}^{(p_1)}  - c_{\beta, 1}^{(p_2)} \theta_{\beta, 1}^{m_2-m_1}},\\
        k_1&=n_1,\ k_2=-m_1,\ k_3=1,\\
        A_1'&=2\log{\theta_{\alpha, 1}},\ A_2'=2\log{\theta_{\beta, 1}}.
    \end{align*}
\noindent Applying \eqref{height_of_sum}, \eqref{height_of_product}, and \eqref{height_of_exponents}, we get 

\begin{align*}
 h \left(\delta_3  \right)
 &\leq h(c_{\alpha, 1}^{(j_1)}) + h(c_{\alpha, 1}^{(j_2)}) 
+ h(c_{\beta, 1}^{(p_1)})+ h(c_{\beta, 1}^{(p_2)}) + 2\log 2 + (m_1 -m_2) \frac{\log \theta_{\beta, 1}}{2} \\ & \ \ \  + (n_1 -n_2) \frac{\log \theta_{\alpha, 1}}{2}.
\end{align*}


\noindent As noted above, 
\begin{align*}
    \max (n_{1}-n_{2}, m_{1}-m_{2}) \leq 1.54 \times10^{29}  (\log \theta_{\beta,1})^{4} (\log n_1)^2,\\
    \text{ and } \ \ \ \ \ \ h(c_{\alpha, 1}^{(j_1)}) + h(c_{\alpha, 1}^{(j_2)}) + h(c_{\beta, 1}^{(p_1)})+ h(c_{\beta, 1}^{(p_2)}) \leq 10 \log \theta_{\beta, 1}.
\end{align*}  

\noindent This gives
\begin{align*}
 h \left( \delta_3  \right) &\leq 1.55 \times10^{29}  (\log \theta_{\beta,1})^{5} (\log n_1)^2
\end{align*}
and therefore, we may take $A_{3}' = 6.2 \times 10^{29}  (\log \theta_{\beta,1})^{5} (\log n_1)^2$.\\
\noindent Using equation \eqref{m1Rn1_eqn2}, we can take $B=2n_1$. Then 
    \[
  \log|\Gamma_{9}-1|>- 1.3558 \times 10^{43} (1+\log 2n_1) \log \theta_{\alpha, 1} (\log \theta_{\beta,1})^{6} (\log n_1)^{2}.  
    \]
Comparing with \eqref{lfl9_rhs}, we obtain
\[ m_1 < 2.72 \times 10^{43} \log \theta_{\alpha,1}(\log \theta_{\beta,1})^{5} (\log n_1)^{3}.\] 
Now, using equation \eqref{n1Rm1_eqn1}, we get
\[n_1 < 2.73 \times 10^{43} (\log \theta_{\beta,1})^{6} (\log n_1)^{3}.
\]

\noindent On combining all four cases we have $n_1 <  ( 3.02 \times 10^{14} (\log \theta_{\beta,1})^2 \log n_1)^{3}.$
    We now use Lemma \ref{PW} with $a=0, c=3, g= 2.73 \times 10^{43} (\log \theta_{\beta,1})^{6}$.  We thus obtain
\begin{align*}
   n_{1} &\leq 216 \times 2.73 \times 10^{43} (\log \theta_{\beta,1})^{6} (\log ( 3 \times 3.02 \times 10^{14} (\log \theta_{\beta,1})^2 ))^{3}\\
   &\leq 5.897 \times 10^{45} (\log \theta_{\beta,1})^{6} (\log ( 9.06 \times 10^{14} (\log \theta_{\beta,1})^2 ))^{3}.
\end{align*}
Hence, \[ \max (N_{1}, M_{1}) \leq 1.1793 \times 10^{46} \times (\log \theta_{\beta,1})^{6} \times (\log (9.06 \times 10^{14} \times (\log \theta_{\beta,1})^{2} ))^{3}.\]

\section{Proof of Theorem \ref{thm explicit_C}}\label{Sec_thm explicit_C}

By a brute force computer search for solutions where $0 \leq N_1 \leq 400$ and $0 \leq M_1 \leq 340$, we find all solutions and list them in Theorem \ref{thm explicit_C}. Therefore, from now onwards we assume that $N_1 \geq 401$\  (and hence $n_1 \geq 200)$.\\

\noindent For $\alpha=[0;\overline{1,2}]$ and $\beta=[0;\overline{1,3}]$, using Theorem \ref{thm on bound} and equation \eqref{m1Rn1_eqn2}, we get 
\[ m_1 < n_1  \leq 3.9 \times 10^{51}.\]

\noindent Now we consider the inequality \eqref{lfl7_rhs} and recall that 
$\Gamma_{7} = \dfrac{c_{\alpha, 1}^{(j_1)}\theta_{\alpha, 1}^{n_1}}{c_{\beta, 1}^{(p_1)} \theta_{\beta, 1}^{m_1}}.$ Taking $\Lambda := \log \Gamma_{7}$, we get
\[ \Lambda = n_{1} \log \theta_{\alpha, 1} - m_{1} \log \theta_{\beta, 1} + \log \left( \frac{c_{\alpha, 1}^{(j_1)}}{c_{\beta, 1}^{(p_1)}}\right)  
.\]
Assuming that $\left| e^{\Lambda} -1 \right| < \frac{1}{4}$, we get $|\Lambda| < \frac{1}{2}$. Since the inequality $|x| < 2|e^{x}-1|$ is true for all nonzero $x \in (-\frac{1}{2}, \frac{1}{2})$, we get
\[|\Lambda| < 20.24 \theta_{\alpha, 1} \times \max \left( \frac{1}{\theta_{\alpha,1}^{n_{1}-n_{2}}}, \frac{1}{\theta_{\beta, 1}^{m_{1}-m_{2}}} \right).\] 

\noindent If $\Lambda > 0$, then we have the following inequality
\[0 <  n_{1} \frac{\log \theta_{\alpha, 1}}{\log \theta_{\beta, 1}} - m_{1}  + \frac{1}{\log \theta_{\beta, 1}} \log \left( \frac{c_{\alpha, 1}^{(j_1)}}{c_{\beta, 1}^{(p_1)}}\right) < \dfrac{20.24  \theta_{\alpha, 1}}{\log \theta_{\beta, 1}} \times \max \left( \frac{1}{\theta_{\alpha,1}^{n_{1}-n_{2}}}, \frac{1}{\theta_{\beta, 1}^{m_{1}-m_{2}}} \right).\]




\noindent If $\Lambda < 0$, then we have the following inequality
$$ 0 <  m_{1} \frac{\log \theta_{\beta, 1}}{\log \theta_{\alpha, 1}} - n_{1}  + \frac{1}{\log \theta_{\alpha, 1}} \log \left( \frac{c_{\beta, 1}^{(p_1)}}{c_{\alpha, 1}^{(j_1)}}\right) < \dfrac{20.24 \theta_{\alpha, 1}} {\log \theta_{\alpha, 1}} \times \max \left( \frac{1}{\theta_{\alpha,1}^{n_{1}-n_{2}}}, \frac{1}{\theta_{\beta, 1}^{m_{1}-m_{2}}} \right).$$ 

\noindent In particular, for $\Lambda > 0 $, $c_{\alpha, 1}^{(j_1)} = \frac{1+\sqrt{3}}{2\sqrt{3}} $ and $ c_{\beta, 1}^{(p_1)}  = \frac{3+\sqrt{21}}{2\sqrt{21}}$,  we apply Lemma \ref{BD_reduction} with 
$\gamma = \dfrac{\log \theta_{\alpha, 1}}{\log \theta_{\beta, 1}} = \frac{\log (2 + \sqrt 3)}{\log ( 5/2 + \sqrt{21}/2)}$, $\kappa = \dfrac{1}{\log \theta_{\beta, 1}} \log \left( \frac{c_{\alpha, 1}^{(j_1)}}{c_{\beta, 1}^{(p_1)}}\right)$, $A_{1} =\dfrac{20.24 \theta_{\alpha, 1}}{\log \theta_{\beta, 1}}$,  $(A_{2}, k) = (\theta_{\alpha,1}, n_{1}-{n_2}) \text{ or } (\theta_{\beta,1}, m_{1} - m_{2}) $, and $M = 3.9 \times 10^{51}$. The continued fraction expansion of $\gamma$ is $[0; 1, 5, 3, 1, 2, 4, 1, 29, 3, 1, \ldots]$. We consider the $105$\emph{th} convergent of the irrational $\gamma$ such that $q=q_{105} > 6M$. This gives $\varepsilon > 0.35261$ and therefore, either 
$n_{1} - n_{2} \leq \frac{\log(\frac{A_{1} q}{0.35261})}{\log \theta_{\alpha, 1}} \leq 95,$ or $m_{1} - m_{2} \leq \frac{\log(\frac{A_{1} q}{0.35261})}{\log \theta_{\beta, 1}} \leq 80$. For other choices of $\Lambda$, $c_{\alpha, 1}^{(j_1)},$ and $ c_{\beta, 1}^{(p_1)}$, we summarize the results in the table below.
\begin{table}[h!]
\centering
\begin{tabular}{|c|c|c|c|c|c|c|}
\hline
        $c_{\alpha, 1}^{(j_1)}$ & $ c_{\beta, 1}^{(p_1)} $ & $\Lambda$ & $ q $ & $\varepsilon$ & $n_1 - n_2 \leq $ & $m_1 - m_2 \leq $  \\
\hline
\multirow{2}{*}{$\frac{1+\sqrt{3}}{2\sqrt{3}}$} & \multirow{2}{*}{$\frac{3+\sqrt{21}}{2\sqrt{21}}$} & $>0$ & $q_{105}$ & $0.35261$& $n_1 - n_2 \leq 95 $ &  $m_1 - m_2 \leq 80$ \\ \cline{3-7} 
&  & $<0$ & $q_{104}$ & $0.34567$& $n_1 - n_2 \leq 95 $ &  $m_1 - m_2 \leq 80$ \\ \hline
\multirow{2}{*}{$\frac{2+\sqrt{3}}{2\sqrt{3}}$} & \multirow{2}{*}{$\frac{3+\sqrt{21}}{2\sqrt{21}}$} & $>0$ & $q_{105}$ & $0.08583$& $n_1 - n_2 \leq 96 $ &  $m_1 - m_2 \leq 81$ \\ \cline{3-7}
&  & $<0$ & $q_{104}$ & $0.07888$& $n_1 - n_2 \leq 96 $ &  $m_1 - m_2 \leq 81$ \\ \hline
\multirow{2}{*}{$\frac{2+\sqrt{3}}{2\sqrt{3}}$} & \multirow{2}{*}{$\frac{5+\sqrt{21}}{2\sqrt{21}}$} & $>0$ & $q_{105}$ & $0.41249$& $n_1 - n_2 \leq 95 $ &  $m_1 - m_2 \leq 80$ \\ \cline{3-7}
&  & $<0$ & $q_{104}$ & $0.40555$& $n_1 - n_2 \leq 95 $ &  $m_1 - m_2 \leq 80$ \\ \hline
\multirow{2}{*}{$\frac{1+\sqrt{3}}{2\sqrt{3}}$} & \multirow{2}{*}{$\frac{5+\sqrt{21}}{2\sqrt{21}}$} & $>0$ & $q_{105}$ & $0.24746$ & $n_1 - n_2 \leq 95 $ &  $m_1 - m_2 \leq 80$ \\ \cline{3-7}
 &   & $<0$ & $q_{104}$ & $0.24051$ & $n_1 - n_2 \leq 95 $ &  $m_1 - m_2 \leq 80$ \\ \hline
\end{tabular}
\caption{Table 1: Computational results for Theorem \ref{thm explicit_C}. }
\end{table}

\noindent If $|e^{\Lambda} -1 | \geq \frac{1}{4}$, then \[ \frac{1}{4}\leq |e^{\Lambda} -1| < 10.12 \theta_{\alpha,1} \times \max \left( \frac{1}{\theta_{\alpha,1}^{n_{1}-n_{2}}}, \frac{1}{\theta_{\beta, 1}^{m_{1}-m_{2}}} \right).\] Therefore, $\min(m_1 - m_2, n_1 -n_2 ) \leq 3$. Combining this with the results obtained in Table 1, we infer that either $n_{1} - n_{2} \leq 96$ or $m_1 - m_2 \leq 81$.
\noindent Now, when $n_1 - n_2 \leq 96,$ we consider the inequality \eqref{lfl8_rhs} and recall that 
$\Gamma_{8} = \dfrac{c_{\alpha, 1}^{(j_1)}- c_{\alpha, 1}^{(j_2)}\theta_{\alpha, 1}^{n_2 - n_1}}{c_{\beta, 1}^{(p_1)}} \theta_{\alpha, 1}^{n_1} \theta_{\beta, 1}^{-m_1}.$ Taking $\Lambda_{1} := \log \Gamma_{8}$, we get
\[ \Lambda_{1} = n_{1} \log \theta_{\alpha, 1} - m_{1} \log \theta_{\beta, 1} + \log \left( \frac{c_{\alpha, 1}^{(j_1)}- c_{\alpha, 1}^{(j_2)}\theta_{\alpha, 1}^{n_2 - n_1}}{c_{\beta, 1}^{(p_1)}}\right).\]
Assuming that $\left| e^{\Lambda_{1}} -1 \right| < \frac{1}{4}$, we get $|\Lambda_{1}| < \frac{1}{2}$. Since the inequality $|x| < 2|e^{x}-1|$ is true for all nonzero $x \in (-\frac{1}{2}, \frac{1}{2})$, we get
\[|\Lambda_{1}| < \frac{4.42}{\theta_{\beta, 1}^{m_1-m_2}}.\]

\noindent If $\Lambda_{1} > 0$, then we have the following inequality
$$ 0 <  n_{1} \frac{\log \theta_{\alpha, 1}}{\log \theta_{\beta, 1}} - m_{1}  + \frac{1}{\log \theta_{\beta, 1}} \log \left( \frac{c_{\alpha, 1}^{(j_1)}- c_{\alpha, 1}^{(j_2)}\theta_{\alpha, 1}^{n_2 - n_1}}{c_{\beta, 1}^{(p_1)}}\right)  < \dfrac{4.42}{\log \theta_{\beta, 1}} \times \frac{1}{\theta_{\beta,1}^{m_{1}-m_{2}}}.$$ 

\noindent If $\Lambda_{1} < 0$, then we have the following inequality
$$ 0 <  m_{1} \frac{\log \theta_{\beta, 1}}{\log \theta_{\alpha, 1}} - n_{1}  + \frac{1}{\log \theta_{\alpha, 1}} \log \left( \frac{c_{\beta, 1}^{(p_1)}}{c_{\alpha, 1}^{(j_1)}- c_{\alpha, 1}^{(j_2)}\theta_{\alpha, 1}^{n_2 - n_1}}\right) < \dfrac{4.42}{\log \theta_{\alpha, 1}} \times \frac{1}{\theta_{\beta,1}^{m_{1}-m_{2}}}.$$ 

\noindent In particular, for $\Lambda_{1} > 0 $, $c_{\alpha, 1}^{(j_1)} = \frac{1+\sqrt{3}}{2\sqrt{3}}$, $c_{\alpha, 1}^{(j_2)} = \frac{1+\sqrt{3}}{2\sqrt{3}}$, and  $ c_{\beta, 1}^{(p_1)}  = \frac{3+\sqrt{21}}{2\sqrt{21}}$,  we apply Lemma \ref{BD_reduction} with 
$\gamma = \dfrac{\log \theta_{\alpha, 1}}{\log \theta_{\beta, 1}}$, $\kappa_{t} = \dfrac{1}{\log \theta_{\beta, 1}} \log \left( \frac{c_{\alpha, 1}^{(j_1)}- c_{\alpha, 1}^{(j_2)}\theta_{\alpha, 1}^{-t}}{c_{\beta, 1}^{(p_1)}}\right)$, $A_{1} =\dfrac{4.42}{\log \theta_{\beta, 1}}$,  $A_{2} = \theta_{\beta,1}$, $k = m_{1}-m_{2}$, and $M = 3.9 \times 10^{51}$. We consider the $109$\emph{th} convergent of the irrational $\gamma$ such that $q=q_{109} > 6M$. For each possible value of $n_{1} -n_{2} = t= 1, \ldots, 96$, we get $\varepsilon > 0.00288$ and therefore,  $m_{1} - m_{2} \leq \frac{\log(\frac{A_{1} q}{0.00288})}{\log \theta_{\beta, 1}} \leq 83.$ For other choices of $\Lambda_1$, $c_{\alpha, 1}^{(j_1)}, c_{\alpha, 1}^{(j_2)}$, and $c_{\beta, 1}^{(p_1)}$, we summarize the results in the table below.
\begin{table}[H]
\centering
\begin{tabular}{|c|c|c|c|c|c|c|}
\hline
$c_{\alpha, 1}^{(j_1)}$ & $c_{\alpha, 1}^{(j_2)}$ &$ c_{\beta, 1}^{(p_1)} $ & $\Lambda_{1}$ & $ q  $ & $\varepsilon$ & $m_1 - m_2 \leq  $  \\
\hline
\multirow{2}{*}{$\frac{1+\sqrt{3}}{2\sqrt{3}}$} & \multirow{2}{*}{$\frac{1+\sqrt{3}}{2\sqrt{3}}$} & \multirow{2}{*}{$\frac{3+\sqrt{21}}{2\sqrt{21}}$} & $>0$ & $q_{109}$ & $0.00288$& $m_1 - m_2\leq 83 $  \\ \cline{4-7} 
&  &  & $<0$ & $q_{109}$ & $0.00730$& $m_1 - m_2 \leq 83 $ \\ \hline
\multirow{2}{*}{$\frac{1+\sqrt{3}}{2\sqrt{3}}$} & \multirow{2}{*}{$\frac{1+\sqrt{3}}{2\sqrt{3}}$} & \multirow{2}{*}{$\frac{5+\sqrt{21}}{2\sqrt{21}}$} & $>0$ & $q_{109}$ & $0.01266$& $m_1 - m_2 \leq 82 $  \\ \cline{4-7} 
&  &  & $<0$ & $q_{109}$ & $0.00293$& $m_1 - m_2 \leq 84 $ \\ \hline
\multirow{2}{*}{$\frac{1+\sqrt{3}}{2\sqrt{3}}$} & \multirow{2}{*}{$\frac{2+\sqrt{3}}{2\sqrt{3}}$} & \multirow{2}{*}{$\frac{3+\sqrt{21}}{2\sqrt{21}}$} & $>0$ & $q_{109}$ & $0.00719$& $m_1 - m_2 \leq 82 $  \\ \cline{4-7} 
&  &  & $<0$ & $q_{109}$ & $0.00153$& $m_1 - m_2 \leq 84 $ \\ \hline
\multirow{2}{*}{$\frac{1+\sqrt{3}}{2\sqrt{3}}$} & \multirow{2}{*}{$\frac{2+\sqrt{3}}{2\sqrt{3}}$} & \multirow{2}{*}{$\frac{5+\sqrt{21}}{2\sqrt{21}}$} & $>0$ & $q_{109}$ & $0.00362$& $m_1 - m_2 \leq 83 $  \\ \cline{4-7} 
&  &  & $<0$ & $q_{109}$ & $0.00758$& $m_1 - m_2 \leq 83 $ \\ \hline
\multirow{2}{*}{$\frac{2+\sqrt{3}}{2\sqrt{3}}$} & \multirow{2}{*}{$\frac{1+\sqrt{3}}{2\sqrt{3}}$} & \multirow{2}{*}{$\frac{3+\sqrt{21}}{2\sqrt{21}}$} & $>0$ & $q_{110}$ & $0.00043$& $m_1 - m_2 \leq 85 $  \\ \cline{4-7} 
&  &  & $<0$ & $q_{110}$ & $0.00210$& $m_1 - m_2 \leq 84 $ \\ \hline
\multirow{2}{*}{$\frac{2+\sqrt{3}}{2\sqrt{3}}$} & \multirow{2}{*}{$\frac{1+\sqrt{3}}{2\sqrt{3}}$} & \multirow{2}{*}{$\frac{5+\sqrt{21}}{2\sqrt{21}}$} & $>0$ & $q_{109}$ & $0.00725$& $m_1 - m_2 \leq 82 $  \\ \cline{4-7} 
&  &  & $<0$ & $q_{111}$ & $0.00994$& $m_1 - m_2 \leq 84 $ \\ \hline
\multirow{2}{*}{$\frac{2+\sqrt{3}}{2\sqrt{3}}$} & \multirow{2}{*}{$\frac{2+\sqrt{3}}{2\sqrt{3}}$} & \multirow{2}{*}{$\frac{3+\sqrt{21}}{2\sqrt{21}}$} & $>0$ & $q_{109}$ & $0.00153$& $m_1 - m_2 \leq 83 $  \\ \cline{4-7} 
&  &  & $<0$ & $q_{109}$ & $0.00518$ & $m_1 - m_2 \leq 84 $ \\ \hline
\multirow{2}{*}{$\frac{2+\sqrt{3}}{2\sqrt{3}}$} & \multirow{2}{*}{$\frac{2+\sqrt{3}}{2\sqrt{3}}$} & \multirow{2}{*}{$\frac{5+\sqrt{21}}{2\sqrt{21}}$} & $>0$ & $q_{109}$ & $0.00304$ & $m_1 - m_2 \leq 83 $  \\ \cline{4-7} 
&  &  & $<0$ & $q_{109}$ & $0.00182$ & $m_1 - m_2 \leq 84$ \\ \hline
\end{tabular}
\caption{Table 2.1: Computational results for Theorem \ref{thm explicit_C}.}
\end{table}

 \noindent If $\left| e^{\Lambda_{1}} -1 \right| \geq \frac{1}{4}$, then \[\frac{1}{4} \leq |e^{\Lambda_{1}} -1| < \frac{2.21}{\theta_{\beta, 1}^{m_1-m_2}}.\] Therefore, $m_1 -m_2 \leq 1$.

 \noindent Finally, if $n_1 - n_2 \leq 96$, then $m_1 - m_2 \leq 85$.\\
\noindent Next, if $m_1 - m_2 \leq 81,$ then we proceed in a similar manner and summarize the results in the table below.

\begin{table}[H]
\centering
\begin{tabular}{|c|c|c|c|c|c|c|}
\hline
$c_{\alpha, 1}^{(j_1)}$ & $c_{\beta, 1}^{(p_1)}$ & $ c_{\beta, 1}^{(p_2)} $ & $\Lambda_{2}$ & $ q  $ & $\varepsilon$ & $n_1 - n_2 \leq  $  \\
\hline
\multirow{2}{*}{$\frac{1+\sqrt{3}}{2\sqrt{3}}$} & \multirow{2}{*}{$\frac{3+\sqrt{21}}{2\sqrt{21}}$} & \multirow{2}{*}{$\frac{3+\sqrt{21}}{2\sqrt{21}}$} & $>0$ & $q_{109}$ & $0.00573$& $n_1 - n_2 \leq 99 $  \\ \cline{4-7} 
&  &  & $<0$ & $q_{109}$ & $0.00515$& $n_1 - n_2 \leq 99 $ \\ \hline
\multirow{2}{*}{$\frac{1+\sqrt{3}}{2\sqrt{3}}$} & \multirow{2}{*}{$\frac{5+\sqrt{21}}{2\sqrt{21}}$} & \multirow{2}{*}{$\frac{3+\sqrt{21}}{2\sqrt{21}}$} & $>0$ & $q_{110}$ & $0.02362$& $n_1 - n_2 \leq 99 $  \\ \cline{4-7} 
&  &  & $<0$ & $q_{111}$ & $0.02374$& $n_1 - n_2 \leq 99 $ \\ \hline
\multirow{2}{*}{$\frac{1+\sqrt{3}}{2\sqrt{3}}$} & \multirow{2}{*}{$\frac{3+\sqrt{21}}{2\sqrt{21}}$} & \multirow{2}{*}{$\frac{5+\sqrt{21}}{2\sqrt{21}}$} & $>0$ & $q_{109}$ & $0.00093$& $n_1 - n_2 \leq 101 $  \\ \cline{4-7} 
&  &  & $<0$ & $q_{109}$ & $0.00060$& $n_1 - n_2 \leq 100 $ \\ \hline
\multirow{2}{*}{$\frac{1+\sqrt{3}}{2\sqrt{3}}$} & \multirow{2}{*}{$\frac{5+\sqrt{21}}{2\sqrt{21}}$} & \multirow{2}{*}{$\frac{5+\sqrt{21}}{2\sqrt{21}}$} & $>0$ & $q_{109}$ & $0.00069$& $n_1 - n_2 \leq 101 $  \\ \cline{4-7} 
&  &  & $<0$ & $q_{109}$ & $0.00303$& $n_1 - n_2 \leq 99 $ \\ \hline
\end{tabular}
\end{table}

\begin{table}[H]
\centering
\begin{tabular}{|c|c|c|c|c|c|c|}
\hline
$c_{\alpha, 1}^{(j_1)}$ & $c_{\beta, 1}^{(p_1)}$ & $ c_{\beta, 1}^{(p_2)} $ & $\Lambda_{2}$ & $ q  $ & $\varepsilon$ & $n_1 - n_2 \leq  $  \\
\hline

\multirow{2}{*}{$\frac{2+\sqrt{3}}{2\sqrt{3}}$} & \multirow{2}{*}{$\frac{3+\sqrt{21}}{2\sqrt{21}}$} & \multirow{2}{*}{$\frac{3+\sqrt{21}}{2\sqrt{21}}$} & $>0$ & $q_{109}$ & $0.00118$& $n_1 - n_2 \leq 101 $  \\ \cline{4-7} 
&  &  & $<0$ & $q_{109}$ & $0.00049$& $n_1 - n_2 \leq 100 $ \\ \hline
\multirow{2}{*}{$\frac{2+\sqrt{3}}{2\sqrt{3}}$} & \multirow{2}{*}{$\frac{5+\sqrt{21}}{2\sqrt{21}}$} & \multirow{2}{*}{$\frac{3+\sqrt{21}}{2\sqrt{21}}$} & $>0$ & $q_{111}$ & $0.02121$& $n_1 - n_2 \leq 99 $  \\ \cline{4-7} 
&  &  & $<0$ & $q_{109}$ & $0.00360$& $n_1 - n_2 \leq 99 $ \\ \hline
\multirow{2}{*}{$\frac{2+\sqrt{3}}{2\sqrt{3}}$} & \multirow{2}{*}{$\frac{3+\sqrt{21}}{2\sqrt{21}}$} & \multirow{2}{*}{$\frac{5+\sqrt{21}}{2\sqrt{21}}$} & $>0$ & $q_{109}$ & $0.01133$& $n_1 - n_2 \leq 99 $  \\ \cline{4-7} 
&  &  & $<0$ & $q_{109}$ & $0.02401$& $n_1 - n_2 \leq 97 $ \\ \hline
\multirow{2}{*}{$\frac{2+\sqrt{3}}{2\sqrt{3}}$} & \multirow{2}{*}{$\frac{5+\sqrt{21}}{2\sqrt{21}}$} & \multirow{2}{*}{$\frac{5+\sqrt{21}}{2\sqrt{21}}$} & $>0$ & $q_{109}$ & $0.00226$& $n_1 - n_2 \leq 100 $  \\ \cline{4-7} 
&  &  & $<0$ & $q_{109}$ & $0.00391$& $n_1 - n_2 \leq 99$ \\ \hline
\end{tabular}
\caption{Table 2.2: Computational results for Theorem \ref{thm explicit_C}.}
\end{table}

\noindent Finally, if $m_1 - m_2 \leq 81$, then $n_1 - n_2 \leq 101$. 
Combining the results obtained so far and taking the maximum of the bounds, we get $n_1 - n_2 \leq 101$ and $m_1 - m_2 \leq 85$. Now, we consider the inequality \eqref{lfl9_rhs} and recall that 
$\Gamma_{9}=\dfrac{c_{\alpha, 1}^{(j_1)}- c_{\alpha, 1}^{(j_2)}\theta_{\alpha, 1}^{n_2 - n_1}}{c_{\beta, 1}^{(p_1)}  - c_{\beta, 1}^{(p_2)} \theta_{\beta, 1}^{m_2-m_1}} \theta_{\alpha, 1}^{n_1} \theta_{\beta, 1}^{-m_1}.$ Taking $\Lambda_{3} := \log \Gamma_{9}$, we get
\[ \Lambda_{3} = n_{1} \log \theta_{\alpha, 1} - m_{1} \log \theta_{\beta, 1} + \log \left( \dfrac{c_{\alpha, 1}^{(j_1)}- c_{\alpha, 1}^{(j_2)}\theta_{\alpha, 1}^{n_2 - n_1}}{c_{\beta, 1}^{(p_1)}  - c_{\beta, 1}^{(p_2)} \theta_{\beta, 1}^{m_2-m_1}}\right).\]
Assuming that $\left| e^{\Lambda_{3}} -1 \right| < \frac{1}{4}$, we get $|\Lambda_{3}| < \frac{1}{2}$. Since the inequality $|x| < 2|e^{x}-1|$ is true for all nonzero $x \in (-\frac{1}{2}, \frac{1}{2})$, we get
\[|\Lambda_{3}| < \frac{3.74}{\theta_{\beta, 1}^{m_1}}.\]

\noindent If $\Lambda_{3} > 0$, then we have the following inequality
$$ 0 <  n_{1} \frac{\log \theta_{\alpha, 1}}{\log \theta_{\beta, 1}} - m_{1}  + \frac{1}{\log \theta_{\beta, 1}} \log \left( \dfrac{c_{\alpha, 1}^{(j_1)}- c_{\alpha, 1}^{(j_2)}\theta_{\alpha, 1}^{n_2 - n_1}}{c_{\beta, 1}^{(p_1)}  - c_{\beta, 1}^{(p_2)} \theta_{\beta, 1}^{m_2-m_1}}\right)  < \dfrac{3.74}{\log \theta_{\beta, 1}} \times \frac{1}{\theta_{\beta,1}^{m_{1}}}.$$ 

\noindent If $\Lambda_{3} < 0$, then we have the following inequality
$$ 0 <  m_{1} \frac{\log \theta_{\beta, 1}}{\log \theta_{\alpha, 1}} - n_{1}  + \frac{1}{\log \theta_{\alpha, 1}} \log \left( \dfrac{c_{\beta, 1}^{(p_1)}  - c_{\beta, 1}^{(p_2)} \theta_{\beta, 1}^{m_2-m_1}}{c_{\alpha, 1}^{(j_1)}- c_{\alpha, 1}^{(j_2)}\theta_{\alpha, 1}^{n_2 - n_1}}\right) < \dfrac{3.74}{\log \theta_{\alpha, 1}} \times \frac{1}{\theta_{\beta,1}^{m_{1}}}.$$ 

\noindent In particular, for $\Lambda_{3} > 0 $, $c_{\alpha, 1}^{(j_1)} = \frac{1+\sqrt{3}}{2\sqrt{3}}$, $c_{\alpha, 1}^{(j_2)} = \frac{1+\sqrt{3}}{2\sqrt{3}}$, $ c_{\beta, 1}^{(p_1)}  = \frac{3+\sqrt{21}}{2\sqrt{21}}$, and $ c_{\beta, 1}^{(p_2)}  = \frac{3+\sqrt{21}}{2\sqrt{21}}$,  we apply Lemma \ref{BD_reduction} with 
$\gamma = \dfrac{\log \theta_{\alpha, 1}}{\log \theta_{\beta, 1}}$, \\$\kappa_{t, s} = \dfrac{1}{\log \theta_{\beta, 1}} \log \left( \dfrac{c_{\alpha, 1}^{(j_1)}- c_{\alpha, 1}^{(j_2)}\theta_{\alpha, 1}^{-t}}{c_{\beta, 1}^{(p_1)}  - c_{\beta, 1}^{(p_2)} \theta_{\beta, 1}^{-s}}\right)$, $A_{1} =\dfrac{3.74}{\log \theta_{\beta, 1}}$,  $A_{2} = \theta_{\beta,1}$, $k = m_{1}$, and $M = 3.9 \times 10^{51}$. We consider the $112$\emph{th} convergent of the irrational $\gamma$ such that $q=q_{112} > 6M$. For each possible value of $n_{1} -n_{2} = t= 1, \ldots, 101$, and $m_{1} - m_{2} = s= 1, \ldots, 85$, we get $\varepsilon > 0.00009$ and therefore, $m_{1} \leq \frac{\log(\frac{A_{1} q}{0.00009})}{\log \theta_{\beta, 1}} \leq 87.$ For other choices of $\Lambda_{3}$, $c_{\alpha, 1}^{(j_1)},\  c_{\alpha, 1}^{(j_2)}$, $ c_{\beta, 1}^{(p_1)},$ and $ c_{\beta, 1}^{(p_2)}$, we summarize the results in the table below.

For the values of $c_{\alpha, 1}^{(j_1)}, c_{\alpha, 1}^{(j_2)}, c_{\beta, 1}^{(p_1)},$ and $ c_{\beta, 1}^{(p_2)} $, given in the last row of the table below, the case where $n_1 - n_2 = 2$ and $m_1 - m_2 =2$ is treated separately. This is because the value of $\varepsilon$ turns out to be negative for all convergent denominators greater than $6M$.

\begin{table}[H]
\centering
\begin{tabular}{|c|c|c|c|c|c|c|c|}
\hline
$c_{\alpha, 1}^{(j_1)}$ & $c_{\alpha, 1}^{(j_2)}$ & $c_{\beta, 1}^{(p_1)}$ & $ c_{\beta, 1}^{(p_2)} $ & $\Lambda_{3}$ & $ q $ & $\varepsilon$ & $m_1\leq  $  \\
\hline
\multirow{2}{*}{$\frac{1+\sqrt{3}}{2\sqrt{3}}$} & \multirow{2}{*}{$\frac{1+\sqrt{3}}{2\sqrt{3}}$} & \multirow{2}{*}{$\frac{3+\sqrt{21}}{2\sqrt{21}}$} & \multirow{2}{*}{$\frac{3+\sqrt{21}}{2\sqrt{21}}$} & $>0$ & $q_{112}$ & $0.00009$& $m_1 \leq 87 $  \\ \cline{5-8} 
&  &  &  & $<0$ & $q_{113}$ & $0.00006$ & $m_1  \leq 90 $ \\ \hline
\multirow{2}{*}{$\frac{1+\sqrt{3}}{2\sqrt{3}}$} & \multirow{2}{*}{$\frac{2+\sqrt{3}}{2\sqrt{3}}$} & \multirow{2}{*}{$\frac{3+\sqrt{21}}{2\sqrt{21}}$} & \multirow{2}{*}{$\frac{3+\sqrt{21}}{2\sqrt{21}}$} & $>0$ & $q_{113}$ & $0.00009$& $m_1 \leq 89 $  \\ \cline{5-8} 
&  &  &  & $<0$ & $q_{113}$ & $0.00001$& $m_1  \leq 90 $ \\ \hline
\multirow{2}{*}{$\frac{1+\sqrt{3}}{2\sqrt{3}}$} & \multirow{2}{*}{$\frac{1+\sqrt{3}}{2\sqrt{3}}$} & \multirow{2}{*}{$\frac{5+\sqrt{21}}{2\sqrt{21}}$} & \multirow{2}{*}{$\frac{3+\sqrt{21}}{2\sqrt{21}}$} & $>0$ & $q_{112}$ & $0.00006$& $m_1  \leq 87 $  \\ \cline{5-8} 
&  &  &  & $<0$ & $q_{112}$ & $6 \times 10^{-6}$& $m_1  \leq 90 $ \\ \hline
\multirow{2}{*}{$\frac{1+\sqrt{3}}{2\sqrt{3}}$} & \multirow{2}{*}{$\frac{2+\sqrt{3}}{2\sqrt{3}}$} & \multirow{2}{*}{$\frac{5+\sqrt{21}}{2\sqrt{21}}$} & \multirow{2}{*}{$\frac{3+\sqrt{21}}{2\sqrt{21}}$} & $>0$ & $q_{114}$ & $0.00003$& $m_1  \leq 90 $  \\ \cline{5-8} 
&  &  &  & $<0$ & $q_{113}$ & $0.00003$& $m_1\leq 90 $ \\ \hline
\multirow{2}{*}{$\frac{1+\sqrt{3}}{2\sqrt{3}}$} & \multirow{2}{*}{$\frac{2+\sqrt{3}}{2\sqrt{3}}$} & \multirow{2}{*}{$\frac{3+\sqrt{21}}{2\sqrt{21}}$} & \multirow{2}{*}{$\frac{5+\sqrt{21}}{2\sqrt{21}}$} & $>0$ & $q_{113}$ & $ 5 \times 10^{-6}$& $m_1 \leq 90 $  \\ \cline{5-8} 
&  &  &  &$<0$ & $q_{113}$ & $2 \times 10^{-6}$& $m_1 \leq 91 $ \\ \hline
\multirow{2}{*}{$\frac{1+\sqrt{3}}{2\sqrt{3}}$} & \multirow{2}{*}{$\frac{1+\sqrt{3}}{2\sqrt{3}}$} & \multirow{2}{*}{$\frac{3+\sqrt{21}}{2\sqrt{21}}$} & \multirow{2}{*}{$\frac{5+\sqrt{21}}{2\sqrt{21}}$} & $>0$ & $q_{113}$ & $0.00005$ & $m_1 \leq 89 $  \\ \cline{5-8} 
&  &  &  &$<0$ & $q_{113}$ & $0.00010$& $m_1 \leq 89 $ \\ \hline
\multirow{2}{*}{$\frac{1+\sqrt{3}}{2\sqrt{3}}$} & \multirow{2}{*}{$\frac{2+\sqrt{3}}{2\sqrt{3}}$} & \multirow{2}{*}{$\frac{5+\sqrt{21}}{2\sqrt{21}}$} & \multirow{2}{*}{$\frac{5+\sqrt{21}}{2\sqrt{21}}$} & $>0$ & $q_{113}$ & $0.00005$& $m_1  \leq 89 $  \\ \cline{5-8} 
&  &  &  & $<0$ & $q_{113}$ & $0.00003$ & $m_1 \leq 90 $ \\ \hline
\multirow{2}{*}{$\frac{1+\sqrt{3}}{2\sqrt{3}}$} & \multirow{2}{*}{$\frac{1+\sqrt{3}}{2\sqrt{3}}$} & \multirow{2}{*}{$\frac{5+\sqrt{21}}{2\sqrt{21}}$} & \multirow{2}{*}{$\frac{5+\sqrt{21}}{2\sqrt{21}}$} & $>0$ & $q_{113}$ & $8 \times 10^{-6}$ & $m_1 \leq 90 $  \\ \cline{5-8} 
&  &  &  & $<0$ & $q_{113}$ & $0.00016$ & $m_1 \leq 89 $ \\ \hline
\multirow{2}{*}{$\frac{2+\sqrt{3}}{2\sqrt{3}}$} & \multirow{2}{*}{$\frac{1+\sqrt{3}}{2\sqrt{3}}$} & \multirow{2}{*}{$\frac{3+\sqrt{21}}{2\sqrt{21}}$} & \multirow{2}{*}{$\frac{3+\sqrt{21}}{2\sqrt{21}}$} & $>0$ & $q_{112}$ & $0.00012$& $m_1 \leq 86 $  \\ \cline{5-8} 
&  &  &  & $<0$ & $q_{112}$ & $0.00004$& $m_1 \leq 89 $ \\ \hline
\multirow{2}{*}{$\frac{2+\sqrt{3}}{2\sqrt{3}}$} & \multirow{2}{*}{$\frac{2+\sqrt{3}}{2\sqrt{3}}$} & \multirow{2}{*}{$\frac{3+\sqrt{21}}{2\sqrt{21}}$} & \multirow{2}{*}{$\frac{3+\sqrt{21}}{2\sqrt{21}}$} & $>0$ & $q_{113}$ & $0.00014$& $m_1  \leq 88 $  \\ \cline{5-8} 
&  &  &  & $<0$ & $q_{113}$ & $0.00002$& $m_1 \leq 90 $ \\ \hline
\multirow{2}{*}{$\frac{2+\sqrt{3}}{2\sqrt{3}}$} & \multirow{2}{*}{$\frac{1+\sqrt{3}}{2\sqrt{3}}$} & \multirow{2}{*}{$\frac{5+\sqrt{21}}{2\sqrt{21}}$} & \multirow{2}{*}{$\frac{3+\sqrt{21}}{2\sqrt{21}}$} & $>0$ & $q_{113}$ & $0.00006$& $m_1  \leq 89 $  \\ \cline{5-8} 
&  &  &  & $<0$ & $q_{113}$ & $0.00001$& $m_1 \leq 91 $ \\ \hline
\multirow{2}{*}{$\frac{2+\sqrt{3}}{2\sqrt{3}}$} & \multirow{2}{*}{$\frac{2+\sqrt{3}}{2\sqrt{3}}$} & \multirow{2}{*}{$\frac{5+\sqrt{21}}{2\sqrt{21}}$} & \multirow{2}{*}{$\frac{3+\sqrt{21}}{2\sqrt{21}}$} & $>0$ & $q_{113}$ & $0.00001$& $m_1  \leq 90 $  \\ \cline{5-8} 
&  &  &  & $<0$ & $q_{113}$ & $0.00001$& $m_1  \leq 90$ \\ \hline
\multirow{2}{*}{$\frac{2+\sqrt{3}}{2\sqrt{3}}$} & \multirow{2}{*}{$\frac{2+\sqrt{3}}{2\sqrt{3}}$} & \multirow{2}{*}{$\frac{3+\sqrt{21}}{2\sqrt{21}}$} & \multirow{2}{*}{$\frac{5+\sqrt{21}}{2\sqrt{21}}$} & $>0$ & $q_{113}$ & $0.00023$& $m_1  \leq 88 $  \\ \cline{5-8} 
&  &  &  & $<0$ & $q_{113}$ & $0.00018$ & $m_1  \leq 89 $ \\ \hline

\multirow{2}{*}{$\frac{2+\sqrt{3}}{2\sqrt{3}}$} & \multirow{2}{*}{$\frac{1+\sqrt{3}}{2\sqrt{3}}$} & \multirow{2}{*}{$\frac{3+\sqrt{21}}{2\sqrt{21}}$} & \multirow{2}{*}{$\frac{5+\sqrt{21}}{2\sqrt{21}}$} & $>0$ & $q_{113}$ & $2 \times 10^{-6}$ & $m_1 \leq 91 $  \\ \cline{5-8} 
&  &  &  & $<0$ & $q_{113}$ & $0.00003$ & $m_1  \leq 90 $ \\ \hline
\multirow{2}{*}{$\frac{2+\sqrt{3}}{2\sqrt{3}}$} & \multirow{2}{*}{$\frac{1+\sqrt{3}}{2\sqrt{3}}$} & \multirow{2}{*}{$\frac{5+\sqrt{21}}{2\sqrt{21}}$} & \multirow{2}{*}{$\frac{5+\sqrt{21}}{2\sqrt{21}}$} & $>0$ & $q_{113}$ & $0.00001$ & $m_1  \leq 90 $  \\ \cline{5-8} 
&  &  &  & $<0$ & $q_{113}$ & $9 \times 10^{-6}$ & $m_1 \leq 91$ \\ \hline
\multirow{2}{*}{$\frac{2+\sqrt{3}}{2\sqrt{3}}$} & \multirow{2}{*}{$\frac{2+\sqrt{3}}{2\sqrt{3}}$} & \multirow{2}{*}{$\frac{5+\sqrt{21}}{2\sqrt{21}}$} & \multirow{2}{*}{$\frac{5+\sqrt{21}}{2\sqrt{21}}$} & $>0$ & $q_{114}$ & $0.00002$ & $m_1 \leq 90 $  \\ \cline{5-8} 
&  &  &  & $<0$ & $q_{114}$ & $0.00004$ & $m_1  \leq 90$ \\ \hline
\end{tabular}
\caption{Table 3: Computational results for Theorem \ref{thm explicit_C}.}
\end{table}

\noindent When $n_1 - n_2 = 2$ and $m_1 - m_2 =2$, we consider the inequality \eqref{3rd Ineq} and obtain the following
\[ \left|\theta_{\beta, 1}^{m_1} - \theta_{\alpha, 1}^{n_1}\right| \leq 0.77.\]

\noindent Dividing by $\theta_{\alpha, 1}^{n_1}$, we get
$ \left| \theta_{\beta, 1}^{m_1} \theta_{\alpha, 1}^{- n_1} - 1 \right| \leq \frac{0.77}{\theta_{\alpha, 1}^{n_1}}.$\\
\noindent If $\left| \theta_{\beta, 1}^{m_1} \theta_{\alpha, 1}^{- n_1} -1 \right| < \frac{1}{4}$, then $\left| \dfrac{\log \theta_{\alpha, 1}}{\log \theta_{\beta, 1}} - \dfrac{m_1}{n_1}\right| < \dfrac{1.54}{n_1 \theta_{\alpha, 1}^{n_1} \log \theta_{\beta, 1}}$. For all $n_{1} \in \mathbb{N}$,\\

\noindent we have $\dfrac{1.54}{n_1 \theta_{\alpha, 1}^{n_1} \log \theta_{\beta, 1}} < \dfrac{1}{2n_{1}^{2}},$ therefore,

\[ \left| \dfrac{\log \theta_{\alpha, 1}}{\log \theta_{\beta, 1}} - \dfrac{m_1}{n_1}\right| < \dfrac{1}{2n_{1}^{2}},\] and this implies that $\dfrac{m_1}{n_1}$ has to be a convergent of the irrational $\dfrac{\log \theta_{\alpha, 1}}{\log \theta_{\beta, 1}}$. 

\noindent Recall that $n_{1} \leq 3.9 \times 10^{51}$. Using this, the convergent with denominator less than or equal to $3.9 \times 10^{51}$ closest to $\dfrac{\log \theta_{\alpha, 1}}{\log \theta_{\beta, 1}}$ could either be $\dfrac{p_{101}}{q_{101}}$ or $\dfrac{p_{102}}{q_{102}}$.
Since $\min \left( \left|\dfrac{\log \theta_{\alpha, 1}}{\log \theta_{\beta, 1}} - \dfrac{p_{101}}{q_{101}} \right|, \left|\dfrac{\log \theta_{\alpha, 1}}{\log \theta_{\beta, 1}} - \dfrac{p_{102}}{q_{102}} \right|\right) = \left|\dfrac{\log \theta_{\alpha, 1}}{\log \theta_{\beta, 1}} - \dfrac{p_{102}}{q_{102}} \right|$, we have
$9 \times 10^{-104} \leq \left|\dfrac{\log \theta_{\alpha, 1}}{\log \theta_{\beta, 1}} - \dfrac{p_{102}}{q_{102}} \right| \leq \dfrac{1.54}{n_1 \theta_{\alpha, 1}^{n_1} \log \theta_{\beta, 1}}$.
This implies $n_1 \leq 176.$\\
\noindent If $ \dfrac{1}{4} \leq \left| \theta_{\beta, 1}^{m_1} \theta_{\alpha, 1}^{- n_1} -1 \right| $, then $\dfrac{1}{4} \leq  \dfrac{0.77}{\theta_{\alpha, 1}^{n_1}}$ which gives $n_1 \leq 0.86$.\\
Therefore, when $n_1 - n_2 =2 = m_1 - m_2$, we get $n_1 \leq 176.$\\
\noindent If $|e^{\Lambda_3} - 1| \geq \frac{1}{4}$, then \[ \frac{1}{4} \leq |e^{\Lambda_3} - 1| \leq \frac{1.87}{\theta_{\beta, 1}^{m_1}}.\]
Therefore, $m_{1} \leq 1$. Finally, combining all the bounds on $m_1$, we have $m_1 \leq 91$. Now, using inequality \eqref{n1Rm1_eqn1}, we have $n_1 \leq 109$. 
Combining the bounds on $n_1$, we have $n_1 \leq 176$. Likewise, we proceed in the other three cases (Cases (i), (ii), and (iii)) and get a bound which is lower than $n_1 \leq 176$. Finally, we get $N_1 \leq 353$, which contradicts our assumption that $N_1 \geq 401$. Hence Theorem \ref{thm explicit_C}.

\section{Appendix}\label{appendix}

\begin{enumerate}[(i)]
    \item For $\alpha = [0; \overline{1,2}]$ and $\beta=[0;\overline{1,4}]$, we have
$$\mathcal{C}=\left\{ -2,0,5,6,10,40\right\}.$$
The following list shows two distinct representations of each $c \in \mathcal{C}:$
\begin{align*}       
  q_{\alpha, 3} - q_{\beta, 3} = 4-6 &=\boxed{-2}  = 3-5 =q_{\alpha,2} - q_{\beta,2},\\
 q_{\alpha, 1} - q_{\beta, 1} = 1-1 &=\boxed{0}  = 1-1 =q_{\alpha,0} - q_{\beta,0},\\
 q_{\alpha, 9} - q_{\beta, 7} = 209-204 &=\boxed{5} = 11-6 =q_{\alpha,4} - q_{\beta,3},\\
q_{\alpha, 6} - q_{\beta, 5} = 41-35 &=\boxed{6}  = 11-5 =q_{\alpha,4} - q_{\beta,2},\\
q_{\alpha, 5} - q_{\beta, 2} = 15-5 &=\boxed{10}  = 11-1 =q_{\alpha,4} - q_{\beta,1},\\
q_{\alpha, 9} - q_{\beta, 6} = 209-169 &=\boxed{40}  = 41-1 =q_{\alpha,6} - q_{\beta,1}.
\end{align*}
 \item  For $\alpha = [0; \overline{1,2}]$ and $\beta=[0;\overline{1,5}]$, we have
$$\mathcal{C}=\left\{ -37,-3,0,8\right\}.$$
The following list shows two distinct representations of each $c \in \mathcal{C}:$
\begin{align*}       
  q_{\alpha, 4} - q_{\beta, 5} = 11-48 &=\boxed{-37}  = 4-41 =q_{\alpha,3} - q_{\beta,4},\\
 q_{\alpha, 3} - q_{\beta, 3} = 4-7 &=\boxed{-3}  = 3-6 =q_{\alpha,2} - q_{\beta,2},\\
 q_{\alpha, 1} - q_{\beta, 1} = 1-1 &=\boxed{0} = 1-1 =q_{\alpha,0} - q_{\beta,0},\\
q_{\alpha, 7} - q_{\beta, 5} = 56-48 &=\boxed{8}  = 15-7 =q_{\alpha,5} - q_{\beta,3}.
\end{align*}

 \item For $\alpha = [0; \overline{1,3}]$ and $\beta=[0;\overline{1,4}]$, we have
$$\mathcal{C}=\left\{ -5,-1,0,18,86\right\}.$$
The following list shows two distinct representations of each $c \in \mathcal{C}:$
\begin{align*}       
  q_{\alpha, 5} - q_{\beta, 4} = 24-29 &=\boxed{-5}  = 1-6 =q_{\alpha,0} - q_{\beta,3},\\
 q_{\alpha, 3} - q_{\beta, 3} = 5-6 &=\boxed{-1}  = 4-5 =q_{\alpha,2} - q_{\beta,2},\\
 q_{\alpha, 1} - q_{\beta, 1} = 1-1 &=\boxed{0} = 1-1 =q_{\alpha,0} - q_{\beta,0},\\
q_{\alpha, 5} - q_{\beta, 3} = 24-6 &=\boxed{18}  = 19-1 =q_{\alpha,4} - q_{\beta,0},\\
q_{\alpha, 7} - q_{\beta, 4} = 115-29 &=\boxed{86}  = 91-5 =q_{\alpha,6} - q_{\beta,2}.
\end{align*}
 
 \item For $\alpha = [0; \overline{1,3}]$ and $\beta=[0;\overline{1,5}]$, we have
$$\mathcal{C}=\left\{ -166,-2,0,18\right\}.$$
The following list shows two distinct representations of each $c \in \mathcal{C}:$
\begin{align*}       
  q_{\alpha, 10} - q_{\beta, 9} = 2089-2255 &=\boxed{-166}  = 115-281 =q_{\alpha,7} - q_{\beta,6},\\
 q_{\alpha, 3} - q_{\beta, 3} = 5-7 &=\boxed{-2}  = 4-6 =q_{\alpha,2} - q_{\beta,2},\\
 q_{\alpha, 1} - q_{\beta, 1} = 1-1 &=\boxed{0} = 1-1 =q_{\alpha,0} - q_{\beta,0},\\
q_{\alpha, 5} - q_{\beta, 2} = 24-6 &=\boxed{18}  = 19-1 =q_{\alpha,4} - q_{\beta,0}.
\end{align*}
 
 \item For $\alpha = [0; \overline{1,4}]$ and $\beta=[0;\overline{1,5}]$, we have
$$\mathcal{C}=\left\{ -6,-1,0,28,163\right\}.$$
The following list shows two distinct representations of each $c \in \mathcal{C}:$
\begin{align*}       
  q_{\alpha, 5} - q_{\beta, 4} = 35-41 &=\boxed{-6}  = 1-7 =q_{\alpha,1} - q_{\beta,3},\\
 q_{\alpha, 3} - q_{\beta, 3} = 6-7 &=\boxed{-1}  = 5-6 =q_{\alpha,2} - q_{\beta,2},\\
 q_{\alpha, 1} - q_{\beta, 1} = 1-1 &=\boxed{0} = 1-1 =q_{\alpha,0} - q_{\beta,0},\\
q_{\alpha, 5} - q_{\beta, 3} = 35-7 &=\boxed{28}  = 29-1 =q_{\alpha,4} - q_{\beta,0},\\
q_{\alpha, 7} - q_{\beta, 4} = 204-41 &=\boxed{163}  = 169-6 =q_{\alpha,6} - q_{\beta,2}.
\end{align*}

\end{enumerate}
\section*{Acknowledgement}
The author sincerely thanks his supervisor Dr. Divyum Sharma for her valuable suggestions, feedback and contributions during the making of this paper.

\bibliographystyle{plain}
\bibliography{references}  
\end{document}